\documentclass[times]{oupau}
\usepackage{amsmath}
\usepackage{amsthm}
\usepackage{amsmath,amscd}
\usepackage{graphicx}%for the scalar multiplication of the module
\usepackage{lmodern}
\usepackage{fix-cm}
\usepackage{bbold}            % 载入字体
\DeclareSymbolFont{bbold}{U}{bbold}{m}{n}
\DeclareSymbolFontAlphabet{\mathbbold}{bbold}
\usepackage[utf8]{inputenc}
\usepackage{amssymb}
\usepackage[hidelinks]{hyperref}
\usepackage{cleveref}
\usepackage{enumitem}
\usepackage[all]{xy}
\usepackage{bm}
\usepackage{aliascnt}
\numberwithin{equation}{section}

\theoremstyle{definition}
\newtheorem{rmk}{Remark}[section]
\newtheorem{defi}{Definition}[section]

\newcommand{\R}{\mathbb{R}}
\newcommand{\Z}{\mathbb{Z}}
\newcommand{\Q}{\mathbb{Q}}

\newcommand{\C}{\mathbb{C}}
\newcommand{\F}{\mathbb{F}}

\newcommand{\E}{\mathbb{E}}

\newcommand{\mf}{\mathfrak}

\newcommand{\bbone}{\mathbbold{1}}
\newcommand{\la}{\lambda}
\newcommand{\La}{\Lambda}

\newcommand{\even}{\text{even}}
\newcommand{\odd}{\text{odd}}
\newcommand{\sq}{\text{sq}}

\newcommand{\ot}{\otimes}

\renewcommand{\L}{\Lambda}

\newcommand{\Y}{\mathbb{Y}}

\DeclareMathOperator{\Hom}{Hom}

\DeclareMathOperator{\Aut}{Aut}

\DeclareMathOperator{\Alt}{Alt}
\DeclareMathOperator{\Her}{Her}

\DeclareMathOperator{\diag}{diag}

\DeclareMathOperator{\Mat}{Mat}

\DeclareMathOperator{\alt}{alt}
\DeclareMathOperator{\her}{her}
\DeclareMathOperator{\Nm}{Nm}
\DeclareMathOperator{\cok}{cok}

\newcommand{\Inj}{\operatorname{Inj}}
\newcommand{\Sur}{\operatorname{Sur}}

\newcommand{\GL}{\mathrm{GL}}

\renewcommand{\l}{\lambda}

\usepackage{color}



\begin{document}

% Enter full title and short title for running headers
\title{Groups With Pairings, Hall Modules, and Hall-Littlewood Polynomials}
\shorttitle{Groups With Pairings, Hall Modules, and Hall-Littlewood Polynomials}

% Enter the publication year and the ID number of the paper
\volumeyear{2025}
\paperID{rnn000}

% Author name(s)
\author{Jiahe Shen\affil{1} and Roger Van Peski\affil{1}}
% Abbreviated author name for running headers
\abbrevauthor{Shen, J. and Van Peski, R.}
% Abbreviated author name for first page header
\headabbrevauthor{Shen, J. and Van Peski, R.}

\address{%
\affilnum{1}Department of Mathematics, Columbia University in the City of New York, NY, USA
}

% Address / e-mail address of corresponding author
\correspdetails{js6157@columbia.edu, rv2549@columbia.edu}

% Received/revised/accepted dates will be entered by the publisher during production of an accepted paper. Please do not edit these placeholders for submission.
\received{1 Month 2025}
\revised{11 Month 20XX}
\accepted{21 Month 20XX}

% Enter details of editor communicating this article
\communicated{A. Editor}

\begin{abstract}
We relate the combinatorics of Hall-Littlewood polynomials to that of abelian $p$-groups with alternating or Hermitian perfect pairings. Our main result is an analogue of the classical relationship between the Hall algebra of abelian $p$-groups (without pairings) and Hall-Littlewood polynomials. Specifically, we define a module over the classical Hall algebra with basis indexed by groups with pairings, and explicitly relate its structure constants to Hall-Littlewood polynomials at different values of the parameter $t$.
    
We also show certain expectation formulas with respect to Cohen-Lenstra type measures on groups with pairings. In the alternating case this gives a new and simpler proof of previous results of Delaunay-Jouhet \cite{delaunay2014p}.
\end{abstract}

\maketitle

\section{Introduction}

The Hall-Littlewood symmetric polynomials are classically known to encode information about subgroups and extensions of abelian $p$-groups (and more general modules), most famously by their isomorphism with the so-called Hall algebra \cite[Chapter II]{Macdonald}. In recent works such as Delaunay-Jouhet \cite{delaunay2014p}, Fulman-Kaplan \cite{fulman2018random}, Cuenca-Olshanski \cite{cuenca2023mackey}, and the authors \cite{shen2024non}, it has become apparent that Hall-Littlewood polynomials also appear in similar problems involving groups with additional structure of bilinear pairings, and related structured random matrices. However, a basic theory like the classical Hall algebra/Hall-Littlewood connection has remained lacking. The main purpose of this paper is provide it, which we do in \Cref{thm: Hall and Symmetric functions}. 

We also include a logically independent---but thematically related---result, \Cref{thm: expectation of hom}, which relates Hall-Littlewood polynomials to certain `Hom-moments' of random modules with alternating or Hermitian pairings, or no pairings at all. Only the Hermitian case is new, but the parallelism of formulas with \Cref{thm: Hall and Symmetric functions} seems worth highlighting. Our proof is also simpler than the previous proof of Delaunay-Jouhet \cite{delaunay2014p} in the known cases, and clarifies the origin of those results, to us at least.

Some of our motivation to understand the context for these results came from the simple fact that groups with pairings are basic objects and it seems natural to see how far the elegant classical Hall algebra story extends. Some came also from the appearance of combinatorics of such groups and probability distributions on them in recent works such as Bhargava et al. \cite{bhargava2013modeling}, the aforementioned Delaunay-Jouhet, Lee \cite{lee2022universality}, Nguyen-Wood \cite{nguyen2022local} and others, motivated by number theory.

\subsection{The Hall algebra and modules over it}
Both the classical Hall algebra results and the present ones feature \emph{Hall-Littlewood polynomials}, a family of symmetric polynomials $P_\la(x_1,\ldots,x_n;t)$ in $n$ variables, indexed by integer partitions $\la = (\la_1,\ldots,\la_n)$ with at most $n$ parts $\la_1 \geq \la_2 \geq \ldots \geq \la_n\geq 0$, which feature an extra parameter $t$ which we take to be real. Explicitly, they are given by 
\begin{equation}
    P_\la(x_1,\ldots,x_n;t) := \frac{1}{V_\la(t)} \sum_{\sigma \in S_n} \sigma\left(x_1^{\la_1}\cdots x_n^{\la_n} \prod_{1 \leq i < j \leq n} \frac{x_i-tx_j}{x_i-x_j}\right),
\end{equation}
where $\sigma$ acts by permuting the variables and $V_\la$ is the normalizing constant making the polynomial monic. 

% They form a basis for symmetric polynomials in $n$ variables, meaning that any product of two of them expands as a linear combination
% \begin{equation}
%     P_\la(x_1,\ldots,x_n;t) \cdot P_\mu(x_1,\ldots,x_n;t) = \sum_{\nu} c_{\la,\mu}^\nu(t) P_\nu(x_1,\ldots,x_n;t)
% \end{equation}
% for certain structure constants $c_{\la,\mu}^\nu(t)$, which deform the classical Littlewood-Richardson coefficients at $t=0$. 

The relation to the combinatorics of abelian groups is as follows. Let $p$ be a prime number. Any finite abelian $p$-group is isomorphic to $G_\la := \bigoplus_i \Z/p^{\la_i}\Z$, for some partition $\la$. Defining generators $u_\la$ for each integer partition, one may define a multiplication
\begin{equation}\label{eq:hall_mult_intro}
    u_\la u_\mu = \sum_{\nu \in \Y} G_{\la,\mu}^\nu u_\nu 
\end{equation}
on these generators, where 
\begin{equation}
    G_{\la,\mu}^\nu :=  \#\{H \leq G_\nu: H \cong G_\la, G_\nu/H \cong G_\mu\}
\end{equation}
are subgroup counts and $\Y$ is the set of all partitions (however, only finitely many terms in the sum will be nonzero). A non-immediate fact is that the multiplication \eqref{eq:hall_mult_intro} then yields an associative algebra structure on $\Z[\{u_\la: \la \in \Y\}]$, called the \emph{Hall algebra}. A classical theorem states that there is an algebra isomorphism from the Hall algebra to the ring $\La$ of symmetric functions\footnote{Essentially, symmetric polynomials in infinitely many variables $x_1,x_2,\ldots$, see \Cref{defi: symmetric functions} for details.} 
\begin{equation}\label{eq:hall_alg_iso_intro}
   u_\la \mapsto p^{-n(\la)} P_\la(x_1,x_2,\ldots;t=1/p),
\end{equation}
which sends these generators to normalized Hall-Littlewood polynomials; here $n(\la) = \sum_{i \geq 1} (i-1)\la_i$.

% %Furthermore, a classical theorem \cite[Chapter V]{mac} states that 
% \begin{equation}
%     c_{\la,\mu}^\nu(t=1/p)p^{-n(\nu)+n(\la)+n(\mu)} = \#\{H \leq G_\nu: H \cong G_\la, G_\nu/H \cong G_\mu\}
% \end{equation}
% where $n(\la) := \sum_{i \geq 1} (i-1)\la_i$,  %One can see from this that the subgroup-counts on the right-hand side, usually denoted $G_{\la,\mu}^\nu$, satisfy an associativity property which makes

%In our main result, \Cref{thm: Hall and Symmetric functions}, we develop the analogues of the Hall algebra/Hall-Littlewood isomorphism (recalled in detail in \Cref{cor: Hall and symmetric functions}) in two settings, namely groups equipped with perfect pairings which are (1) alternating or (2) Hermitian. A logically independent---but thematically related---result, \Cref{thm: expectation of hom}, relates Hall-Littlewood polynomials to certain `Hom-moments' of random modules with alternating and Hermitian pairings. Only the Hermitian case is new, but the parallelism of formulas with \Cref{thm: Hall and Symmetric functions} made both cases seem worth highlighting, and our (identical) proof of both is substantially simpler than the previous proof of Delaunay-Jouhet \cite{delaunay2014p} in the alternating case.

Finite abelian $p$-groups are just modules over the $p$-adic integers $\Z_p$, and have the same combinatorics as modules over any non-archimedean local field with finite residue field. Everything stated above is well known to hold in this context (\Cref{cor: Hall and symmetric functions}), so we work in the following settings:

\begin{enumerate}
\item (Alternating case) 
Let 

$F$ be a non-archimedean local field with characteristic zero;%\footnote{The characteristic $0$ assumption is merely technical, see \Cref{rmk:why_char_0} for an explanation of where it is needed in the proofs.};
    
$\mathfrak{o}$ be the ring of integers of $F$;

$\mathfrak{p}$ be the maximal ideal of $\mathfrak{o}$;

$\pi\in \mathfrak{p}$ be a generator of $\mathfrak{p}$;

$k=\mathfrak{o}/\mathfrak{p}$ be the residue field of $F$;

$q=\#k$ be the order of the residue field of $F$;

$|\cdot|:F^\times\rightarrow q^{\Z},x\mapsto q^{-v(x)}$ be the absolute value defined over nonzero elements of $F$; 

%$\Alt_n(F),\Alt_n(\mathfrak{o})$ be the set of alternating\footnote{i.e. $A=-A^T$.} matrices with entries in $F,\mathfrak{o}$, respectively.

\item(Hermitian case) Let 

$F$ be a non-archimedean local field with characteristic zero;

$\mathfrak{o}$ be the ring of integers of $F$;

$\mathfrak{p}$ be the maximal ideal of $\mathfrak{o}$;

$*$ be an involution of $F$, i.e., an isomorphism of order $2$;

$F_0$ be the fixed field of the involution, and suppose furthermore that $F/F_0$ is unramified, so $F_0\cap\mathfrak{o}$ is the ring of integers of $F_0$;

$\pi\in F_0\cap\mathfrak{p}$ be a generator of $\mathfrak{p}$;

$k=\mathfrak{o}/\mathfrak{p}$ be the residue field of $F$;

$q^2=\#k$ be the order of the residue field of $F$;

$|\cdot|:F^\times\rightarrow q^{\Z},x\mapsto q^{-2v(x)}$ be the absolute value defined over nonzero elements of $F$;

%$B^* = *(B^T)$ where $*$ acts entrywise on the matrix;

%$\Her_n(F),\Her_n(\mathfrak{o})$ be the set of Hermitian matrices with entries in $F,\mathfrak{o}$, respectively.

\end{enumerate}

% \begin{example}
%     In the alternating case, for any prime $p$ one may take $F = \Q_p, \mf{o} = \Z_p$ to be the $p$-adic numbers and integers respectively, in which case $\pi = p$ and $k = \F_p$. In the Hermitian case, one may take $F = \Q_p[\sqrt{d}]$ where $d \in \Z \subset \Q_p$ is not a square in $\Z/p\Z$, and $F_0=\Q_p$; then $*$ acts by $*(a+b\sqrt{d}) = a-b\sqrt{d}$. In particular if $p \equiv 3 \pmod{4}$ one may take $F=\Q_p[\sqrt{-1}]$.
% \end{example}

Let $M$ be a finite $\mathfrak{o}$-module, equipped with a bilinear alternating (resp. Hermitian) pairing $\langle\cdot,\cdot\rangle:M\times M\rightarrow F/\mathfrak{o}$. We say $M$ is an \emph{alternating (resp. Hermitian) $\mathfrak{o}$-module} if the pairing $\langle\cdot,\cdot\rangle$ is perfect. Every alternating $\mathfrak{o}$-module $M$ has even length and is isomorphic to one of the form
\begin{equation}
M\cong(\mathfrak{o}/\mathfrak{p}^{\lambda_1}\oplus\mathfrak{o}/\mathfrak{p}^{\lambda_1})\oplus\cdots\oplus(\mathfrak{o}/\mathfrak{p}^{\lambda_n}\oplus\mathfrak{o}/\mathfrak{p}^{\lambda_n})
\end{equation}
for some $\lambda_1\ge\ldots\ge\lambda_n\ge 1$ with the pairing $\langle\cdot,\cdot\rangle$ given by $\langle x,y\rangle=x^T(\pi_\lambda^{\alt})^{-1}y$ for any $x=(x_1,\ldots,x_{2n}),y=(y_1,\ldots,y_{2n})\in M$, where
\begin{equation}
\pi^{\alt}_\lambda=\diag\left({\begin{pmatrix} 0 & \pi^{\lambda_1} \\ -\pi^{\lambda_1} & 0\end{pmatrix}},\ldots,{\begin{pmatrix} 0 & \pi^{\lambda_n} \\ -\pi^{\lambda_n} & 0\end{pmatrix}}\right).
\end{equation}
Here and below, we use $x=(x_1,\ldots,x_{m})$ for elements of $\mf{o}^{m}$ or a quotient thereof. When acted on by a matrix, these should be viewed as column vectors, despite being written as rows; when we wish to refer to a row vector in $\Hom(\mf{o}^{m},F)$ we will write it as $x^T=(x_1,\ldots,x_m)^T$.

Likewise, every Hermitian $\mathfrak{o}$-module $M$ is isomorphic to one of the form
\begin{equation}
M\cong\mathfrak{o}/\mathfrak{p}^{\lambda_1}\oplus\cdots\oplus\mathfrak{o}/\mathfrak{p}^{\lambda_n}
\end{equation}
for some $\lambda_1\ge\ldots\ge\lambda_n\ge 1$ with the pairing $\langle\cdot,\cdot\rangle$ given by $\langle x,y\rangle= x^*\pi_\lambda^{-1}y$ for any $x=(x_1,\ldots,x_n),y=(y_1,\ldots,y_n)\in M$, where $\pi_\lambda=\diag(\pi^{\lambda_1},\ldots,\pi^{\lambda_n})$. Here $x^* = *(x^T)$ where $*$ is the `conjugation' involution defined above.

We say $M$ is an alternating (resp. Hermitian) $\mathfrak{o}$-module \emph{of type} $\l$ for the above case. For any $\mathfrak{o}$-submodule $N$ of $M$, we write
\begin{equation}
    N^\perp := \{x \in M: \langle x,y \rangle = 0 \text{ for all }y \in N\}.
\end{equation}
We call $N$ isotropic if $N \subset N^\perp$, equivalently $\langle x,y \rangle = 0$ for all $x,y \in N$. Denote by $G_{\mu,\nu}^{\alt,\lambda}(\mathfrak{o})$ (resp. $G_{\mu,\nu}^{\her,\lambda}(\mathfrak{o})$) the number of submodules $M'\subset M$ such that $M/M'$ is an $\mathfrak{o}$-module of type $\mu$, $M'^\perp$ is isotropic, and $M'/M'^\perp$ is an alternating (resp. Hermitian) $\mathfrak{o}$-module of type $\nu$. Motivated by the structure of the Hall algebra $H(\mathfrak{o})$ and related structure at the level of the Hecke ring studied previously in \cite{shen2024non}, we define $H^{\alt}=H^{\alt}(\mathfrak{o})$ and $H^{\her}=H^{\her}(\mathfrak{o})$ to be the $H(\mathfrak{o})$-modules with free $\Z$-basis $(u_\l^{\alt})_{\l\in\Y}$ and $(u_\l^{\her})_{\l\in\Y}$ respectively, and with $H(\mf{o})$-actions given by 
\begin{equation}
u_\mu u_\nu^{\alt}=\sum_\lambda G_{\mu,\nu}^{\alt,\lambda}(\mathfrak{o})u_\lambda^{\alt},\quad u_\mu  u_\nu^{\her}=\sum_\lambda G_{\mu,\nu}^{\her,\lambda}(\mathfrak{o})u_\lambda^{\her}.    
\end{equation}
The reader may refer to the detailed discussions following \Cref{def:alt_hall_module} (alternating case) and \Cref{def:her_hall_module} (Hermitian case) for further clarification of the above definition. As an analogue of \eqref{eq:hall_alg_iso_intro}, we construct an isomorphism which connects the alternating (resp. Hermitian) Hall module with the ring $\La_\Q$ of symmetric functions with coefficients in $\Q$ (\Cref{defi: symmetric functions}). We also require the subalgebra $\La_\Q^\sq \subset \La_\Q$, which informally is just the `symmetric functions in infinitely many squared variables $x_1^2,x_2^2,\ldots$', see \Cref{def:Lambda_sq}.

\begin{theorem}\label{thm: Hall and Symmetric functions}
In the above setting, we have:
\begin{enumerate}
\item (Alternating case) Define $\Q$-linear mappings
\begin{align}\label{eq: alt_image of ring}
\begin{split}\phi:H(\mathfrak{o})\otimes_\Z\Q&\rightarrow\Lambda_\Q \\ 
 u_\mu &\mapsto q^{-n(\mu)+|\mu|}P_\mu(x_1,x_1q^{-1},x_2,x_2q^{-1},\ldots;q^{-1})   
\end{split}\end{align}
and
\begin{align}\label{eq: alt_image of module}
\begin{split}
\phi^{\alt}:H^{\alt}(\mathfrak{o})\otimes_\Z\Q&\rightarrow\Lambda_\Q\\
 u_\nu^{\alt}&\mapsto q^{-2n(\nu)}P_\nu(x_1,x_2,\ldots;q^{-2}).  
\end{split}
\end{align}
Then we have 
\begin{equation}\label{eq: alt_Satake isomorphism}
\phi(u_\mu)\phi^{\alt}(u_\nu^{\alt})=\phi^{\alt}(u_\mu*u_\nu^{\alt}),\quad\forall \mu,\nu\in \Y.
\end{equation}
\item (Hermitian case) Define $\Q$-linear mappings
\begin{align}\label{eq: her_image of ring}
\begin{split}\phi:H(\mathfrak{o})\otimes_\Z\Q&\rightarrow\Lambda_\Q^\sq \\ 
 u_\mu &\mapsto q^{-2n(\mu)}P_\mu(x_1^2,x_2^2,\ldots;q^{-2})   
\end{split}\end{align}
and
\begin{align}\label{eq: her_image of module}
\begin{split}
\phi^{\her}:H^{\her}(\mathfrak{o})\otimes_\Z\Q&\rightarrow\Lambda_\Q\\
 u_\nu^{\her}&\mapsto (-1)^{n(\nu)}q^{-n(\nu)}P_\nu(x_1,x_2,\ldots;-q^{-1}).  
\end{split}\end{align}
Then we have
\begin{equation}\label{eq: her_Satake isomorphism}
\phi(u_\mu)\phi^{\her}(u_\nu^{\her})=\phi^{\her}(u_\mu*u_\nu^{\her}),\quad\forall \mu,\nu\in \Y.
\end{equation}
\end{enumerate}
Moreover, the mappings $(\phi,\phi^{\alt})$ (resp. $(\phi,\phi^{\her})$) give a module isomorphism from the $H(\mathfrak{o})\otimes_\Z\Q$-module $H^{\alt}(\mathfrak{o})\otimes_\Z\Q$ (resp. $H^{\her}(\mathfrak{o})\otimes_\Z\Q$) to $\Lambda_\Q$, viewed as a module of itself (resp. as a module of $\Lambda_\Q^\sq$). 
\end{theorem}

\begin{rmk}
Note that we use $\phi$ for different maps in the alternating and Hermitian cases above. Both isomorphisms $\phi$ have a slight difference from the standard one, which we recall later in \eqref{eq: standard form of isomorphism}. 
\end{rmk}

To prove \Cref{thm: Hall and Symmetric functions}, we apply results of Hironaka \cite{Hironaka_Hermitian,Hironaka_Hermitian_and_Symmetric_I} and Hironaka-Sato \cite{Hironaka}, which we became aware of and reformulated in \cite{shen2024non}, to relate Hall-Littlewood polynomials to certain modules relating to action of $\GL_N$ on spaces of alternating and Hermitian matrices. The latter can then be related to our Hall module actions.

To roughly presage the argument, we will show that a precise version of the following diagram is commutative. Here curved arrows denote rings acting on modules, and $\La_n$ the ring of symmetric polynomials in $n$ variables.
\begin{equation}\label{eq:commutative diagrams}
\xymatrix@C=0.5em{
H=H(\mathfrak{o})\ar[d]^{(2)} &\overset{(1)}{\scalebox{2}[1]{\(\curvearrowright\)}} & H^{\alt}=H^{\alt}(\mathfrak{o}) \ar[d]^{(3)}\\
\substack{\text{Hecke} \\ \text{ring}} \ar[d]^{(5)}&\overset{(4)}{\scalebox{2}[1]{\(\curvearrowright\)}} &\substack{\text{Alternating} \\ \text{Hecke module}}\ar[d]^{(6)}\\
\Lambda_n\otimes_\Z\Q&\overset{(7)}{\scalebox{2}[1]{\(\curvearrowright\)}} &\Lambda_n\otimes_\Z\Q
} 
% \hspace{10ex} \xymatrix@C=0.5em{
% H=H(\mathfrak{o})\ar[d]^{(2)} &\overset{(1)}{\scalebox{2}[1]{\(\curvearrowright\)}} & H^{\her}=H^{\her}(\mathfrak{o}) \ar[d]^{(3)}\\
% H(G^+,K) \ar[d]^{(5)}&\overset{(4)}{\scalebox{2}[1]{\(\curvearrowright\)}} & H(G^{+{\her}},K)\ar[d]^{(6)}\\
% \L_n^\sq\otimes_\Z\Q&\overset{(7)}{\scalebox{2}[1]{\(\curvearrowright\)}} &\Lambda_n\otimes_\Z\Q}
\end{equation}
For precise definitions and the exact mappings, see \Cref{thm: alt_commutative diagram}; the Hermitian case is similar, see \Cref{thm: her_commutative diagram}. The map (2) is standard, and the map (5) is a slight variant of the standard one in e.g. \cite[Chapter V]{Macdonald}, see \Cref{thm: left half}; their composition is the isomorphism given in \eqref{eq:hall_alg_iso_intro} composed with a projection map $\La \to \La_n$. In \cite{shen2024non} we introduced the lower two module structures (4) and (7), and showed that the diagram given by the lower square commutes. Our progress in the present work is to introduce the maps (1) and show commutativity of the upper square, which in combination with the commuting of the lower square yields \Cref{thm: Hall and Symmetric functions}.

\begin{rmk}\label{rmk:why_char_0}
    The reason for our characteristic $0$ assumption is that the needed inputs from \cite{Hironaka_Hermitian,Hironaka_Hermitian_and_Symmetric_I,Hironaka} are proven under the same assumption. However, we have been informed of work in preparation by C. Cuenca and G. Olshanski which will extend these results to positive characteristic, at least in the odd case. All arguments we make within this paper do not distinguish between characteristic $0$ and positive characteristic, so once that work appears, the same proofs we give here will apply also in positive odd characteristic.
\end{rmk}

\subsection{Random modules with pairings}\label{subsec: Random modules with pairings} 

In Cohen and Lenstra's original paper on heuristic distributions of class groups \cite{cohen-lenstra}, they defined the $u$-probability for partitions corresponding to (usual) $\mathfrak{o}$-modules, where the probability of a random $\mathfrak{o}$-module $M$ has type $\l$ is given by (here $q$ refers to the order of the residue field of $\mathfrak{o}$)
\begin{equation}
    \label{eq:P}
    \mathbf{P}(u;\lambda)=\frac{q^{-u|\lambda|}}{\#\Aut(M_\lambda)}\prod_{j\ge 1}(1-q^{-u-j}).
\end{equation}
The case $u=0$, known as the Cohen-Lenstra measure, is most commonly studied and conjectured to govern distributions of class groups of quadratic imaginary number fields \cite{cohen-lenstra}. The case $u=1$ conjecturally governs statistics of class groups of real quadratic extensions, see Wood \cite{wood2018cohen}. Inspired by the definition, Delaunay \cite{delaunay} defined an alternating analogue 
\begin{equation}
    \label{eq:P_alt}
    \mathbf{P}^{\alt}(u;\lambda)=\frac{q^{-2u|\lambda|}\#M_\lambda}{\#\Aut^A(M_\lambda)}\prod_{j\ge 1}(1-q^{-2u-2j+1}),
\end{equation}
where $M_\la$ is an alternating $\mf{o}$-module of type $\la$ and $\Aut^A(M_\lambda)$ is the set of module automorphisms which respect the pairing. Conjecturally, these model the distribution of Tate-Shafarevich groups of elliptic curves, see Delaunay \cite{delaunay}, Bhargava et al. \cite{bhargava2013modeling}, Smith \cite{smith20172}.

It is natural to define a similar $u$-probability in the Hermitian case by
\begin{equation}
    \label{eq:P_her}
    \mathbf{P}^{\her}(u;\lambda)=\frac{q^{-u|\lambda|}}{\#\Aut^H(M_\lambda)}\prod_{j\ge 1}(1-(-1)^{j-1}q^{-u-j}),
\end{equation}
where $M_\la$ is a Hermitian $\mf{o}$-module and $\Aut^H(M_\la)$ is the set of module automorphisms which respect the pairing.

These measures also govern limiting distributions of cokernels of large random matrices over $\mf{o}$, and limiting distributions of conjugacy classes of matrices over $\F_q$, see \Cref{rmk:u_meas_and_cokernels} and \Cref{rmk:u_probs_in_finite_field} respectively.

Our next result shows that each of the three cases enjoys very similar moment formulas, the first two of which were shown by Delaunay-Jouhet \cite{delaunay2014p} in a slightly different form. When below we write a Hall-Littlewood polynomial $P_\la(a_1,a_2,\ldots;t)$ with an appropriately decaying sequence of infinitely many real variables $a_1,a_2,\ldots$, this simply denotes $\lim_{n \to \infty} P_\la(a_1,a_2,\ldots,a_n;t)$.
% Meanwhile, \cite[(1.6), Chapter II]{Macdonald} computes the number of automorphisms $\#\Aut(M_\lambda)=q^{|\lambda|+2n(\lambda)}\prod_{i\ge 1}(q^{-1};q^{-1})_{m_i(\lambda)}$. 

\begin{theorem}\label{thm: expectation of hom}
Let $\nu=(\nu_1,\nu_2,\ldots,\nu_n)\in\Y$ be a fixed partition, where $\nu_1\ge\ldots\ge\nu_n\ge 1$. Let 
$$M_\nu\cong(\mathfrak{o}/\mathfrak{p}^{\nu_1})\oplus\cdots\oplus(\mathfrak{o}/\mathfrak{p}^{\nu_n})$$
be a usual $\mathfrak{o}$-module of type $\nu$. For any $\mf{o}$-module $M$, let $\Hom(M,M_\nu)$ denote the set of $\mf{o}$-module homomorphisms from $M$ to $M_\nu$ (even when the modules have pairings, we consider all homomorphisms, not just those which respect the pairing).
\begin{enumerate}
\item (No-pairing case) Let $M$ be a random $\mf{o}$-module with distribution $\mathbf{P}(M\text{ has type }\l) = \mathbf{P}(u;\lambda)$. Let $q=\#(\mathfrak{o}/\mathfrak{p})$ be the order of the residue field. Then 
\begin{equation}
    \E[\#\Hom(M,M_\nu)]= \frac{P_\nu(1,q^{-1},\ldots,q^{-u},q^{-1-u},\ldots;q^{-1})}{P_\nu(1,q^{-1},\ldots;q^{-1})}.
\end{equation}

\item (Alternating case) Let $M$ be a random alternating $\mathfrak{o}$-module with distribution $\mathbf{P}(M\text{ has type }\l)=\mathbf{P}^{\alt}(u;\lambda)$. Then
\begin{equation}
    \E[\#\Hom(M,M_\nu)]=\frac{P_\nu(1,q^{-2},\ldots,q^{1-2u},q^{-1-2u},\ldots;q^{-2})}{P_\nu(1,q^{-2},\ldots;q^{-2})},
\end{equation}
where $q=\#(\mf{o}/\mf{p})$ as before.

\item (Hermitian case) Let $M$ be a random Hermitian $\mathfrak{o}$-module with distribution $\mathbf{P}(M\text{ has type }\l)=\mathbf{P}^{\her}(u;\lambda)$. Then we have
\begin{equation}
    \E[\#\Hom(M,M_\nu)]=\frac{P_{\nu^2}(1,-q^{-1},q^{-2},\ldots,-q^{-u},q^{-1-u},-q^{-2-u},\ldots;-q^{-1})}{P_{\nu^2}(1,-q^{-1},q^{-2},\ldots;-q^{-1})},
\end{equation}
where $\nu^2=(\nu_1,\nu_1,\nu_2,\nu_2,\ldots,\nu_n,\nu_n)$ is the partition that consists of two copies of $\nu$. Recall here that $q^2$ is the size of the residue field of $F$, while $q$ is the size of the residue field of $F_0$.
\end{enumerate}
\end{theorem}

Note that the Hall-Littlewood parameters $q^{-1},q^{-2}$, and $-q^{-1}$ are the same as those appearing in the classical Hall algebra isomorphism \eqref{eq:hall_alg_iso_intro} and in the two cases of \Cref{thm: Hall and Symmetric functions}, respectively.  Despite this, our proof of \Cref{thm: expectation of hom} is logically independent of \Cref{thm: Hall and Symmetric functions}. In particular, it is valid in arbitrary characteristic, as the number of $\mf{o}$-module homomorphisms between two modules is well known to depend only on the size of the residue field, see e.g. \cite[Chapter II]{Macdonald}. Our proof is also logically independent of the argument of \cite{delaunay2014p}, which relied on fairly intricate combinatorics. Instead, it comes directly from a certain identity of Hall-Littlewood polynomials. The dual Hall-Littlewood polynomial $Q_\la$ below is a certain constant multiple of $P_\la$, see \Cref{sec:prelim}, and we use the standard notation $\la_i' := \#\{j: \la_j \geq i\}$ for conjugate parts of a partition.

\begin{proposition}\label[proposition]{thm: expectation over CL measure}
Let $t \in (-1,1) \setminus \{0\}$, and $a_1,a_2,\ldots$ be real numbers such that the sum
\begin{equation}\label{eq:hl_measure_intro}
    \sum_{\mu \in \Y} P_\mu(a_1,a_2,\ldots;t)Q_\mu(1,t,t^2,\ldots;t)
\end{equation}
converges and the summands are all nonnegative. Let $\la$ be a random integer partition with $\mathbf{P}(\la = \mu)$ proportional to the $\mu$ summand in \eqref{eq:hl_measure_intro}, for each $\mu \in \Y$. Then for any fixed partition $\nu$,
\begin{equation}\label{eq: expectation over CL measure}
\E[t^{-\sum_i\lambda_i'\nu_i'}]= \frac{P_\nu(a_1,a_2,\ldots,t,t^2,\ldots;t)}{P_\nu(t,t^2,\ldots;t)}. 
\end{equation}
\end{proposition}

Interestingly, the proof of \Cref{thm: expectation over CL measure} is quite short, but the key step (\Cref{lem: sum of skew}) is proven by relating back to combinatorics of $\mf{o}$-modules. Once one reinterprets the various $u$-measures as special cases of the probability measure in \eqref{eq:hl_measure_intro}, the left hand side of \eqref{eq: expectation over CL measure} naturally gives homomorphism counts and one obtains \Cref{thm: expectation of hom}.

\begin{rmk}
    Our interest in \Cref{thm: expectation over CL measure} and its consequences came partly from the fact that such `$t$-moment' expectations as on the left hand side of \eqref{eq: expectation over CL measure}, for a slightly more general family of so-called Hall-Littlewood measures than \eqref{eq:hl_measure_intro}, are a key tool in integrable probability---see Borodin-Corwin \cite{borodin2014macdonald} or Dimitrov \cite{dimitrov2018kpz}. The formulas for these moments established in \cite{borodin2014macdonald} and many works since take the form of contour integrals, and are proven using Macdonald difference operators. Somewhat surprisingly, \Cref{thm: expectation over CL measure} and its proof are fairly general but feature neither. 
\end{rmk}

\begin{rmk}\label{rmk:why_not_symmetric}
    One may also ask whether analogous results to ours exist for modules with symmetric pairings. This may be a more difficult question, for several reasons. Firstly, isomorphism classes of such modules are not simply parametrized by partitions, as the pairings are only equivalent to $\diag(u_1\pi^{-\la_1},\ldots,u_n \pi^{-\la_n})$ where $u_1,\ldots,u_n$ are either $1$ or a quadratic non-residue. This makes the connection to Hall-Littlewood polynomials less clear. For \Cref{thm: Hall and Symmetric functions} we lack the needed input: while some results are proven for symmetric matrices in \cite{Hironaka_Hermitian_and_Symmetric_I}, the formulas are less explicit. 
    
    For analogues of \Cref{thm: expectation of hom}, there is a connection between modules with symmetric pairings and Hall-Littlewood polynomials proven in \cite{fulman2016hall}. However, these measures are not of the type in \Cref{thm: expectation over CL measure}, which come from the Cauchy identity for Hall-Littlewood polynomials, but instead of a form coming from the Littlewood identities.
\end{rmk}

\subsection{Plan of paper} In \Cref{sec:prelim} we give preliminaries on symmetric functions and the classical material on Hall algebras, Hall-Littlewood polynomials, and the Hecke ring (also known as the type $A$ spherical Hecke algebra). We prove the alternating and Hermitian cases of \Cref{thm: Hall and Symmetric functions} in \Cref{sec:alternating} and \Cref{sec:hermitian} respectively, both of which are quite similar. In \Cref{sec:u} we prove \Cref{thm: expectation over CL measure} and use it to deduce \Cref{thm: expectation of hom}. \hyperref[sec:appendix]{Appendix A} gives an alternate direct (but longer) combinatorial proof of the key lemma used to establish \Cref{thm: expectation over CL measure}.

\section{Preliminaries}\label{sec:prelim}

In this section, we give the classical definitions of the Hall algebra $H=H(\mathfrak{o})$, Hecke ring, and Hall-Littlewood polynomials, and state the relations between them. This may all be found in \cite{Macdonald} and we largely follow the notation there. % according to \cite[Chapter II]{Macdonald},  Many of the objects we shall consider in this section will turn out to be parametrized by partitions, which are defined as follows.

\begin{defi}
Let $\Y$ denote the set of \emph{partitions} $\lambda=(\lambda_1,\lambda_2,\ldots)$, which are (finite or infinite) sequences of nonnegative integers $\lambda_1\ge\lambda_2\ge\cdots$ which are eventually zero. We will not distinguish between two such sequences that differ only by a string of zeros at the end. We refer to the integers $\l_i>0$ as the \emph{parts} of $\l$, and set $|\l| := \sum_{i\ge 1}\l_i,n(\lambda):=\sum_{i\ge 1} (i-1)\lambda_i$, and $m_k(\l) = \#\{i\mid \l_i = k\}$. We also write $l=l(\lambda)=\#\{i\mid \l_i >0\}$ as the \emph{length} of $\lambda$, i.e., the number of nonzero parts, and denote the set of partitions with length $\le n$ by $\Y_n$. For $\l,\mu\in\Y$, write $\mu\subset\l$ if $\l_i \geq \mu_i$ for all $i\ge 1$. 
\end{defi}

\begin{defi}
Let $\l\in\Y$ be a partition. The \emph{conjugate} of the partition $\l$ is the partition $\l'=(\l_1',\l_2',\ldots)$ such that $\l_k'=\#\{i\mid\l_i\ge k\}$ for all $k\ge 1$. In particular, $\l_1'=l(\l),\l_1=l(\l'),(\l')'=\l$, and $m_i(\l)=\l_i'-\l_{i+1}'$.
\end{defi}

\subsection{Hall algebra} Fix a non-archimedean local field $F$ with $\mathfrak{o}$ its ring of integers, $\pi$ a generator of the maximal ideal $\mathfrak{p}$, and $q$ the order of the residue field. Any nonzero element $x\in F^\times$ may be written as $x=\pi^my$ with $m\in\Z$ and $y\in\mathfrak{o}^\times$. Define $|\cdot|: F \to \R_{\ge 0}$ by setting $|x| = q^{-m}$ for $x$ as before, and $|0|=0$. Then $|\cdot|$ defines a norm on $F$ and $d(y_1,y_2) :=|y_1-y_2|$ defines a metric. We additionally define $v(x)=m$ for $x$ as above and $v(0)=\infty$, so $|x|=q^{-v(x)}$. 

The definition of the Hall algebra $H=H(\mf{o})$ can be found in \cite[Chapter II]{Macdonald}. Let us restate it here from the beginning. Recall the classification of finite $\mathfrak{o}$-modules under isomorphism:

\begin{proposition}
Every finite $\mathfrak{o}$-module $M$ is isomorphic to one of the form
$$M\cong(\mathfrak{o}/\mathfrak{p}^{\lambda_1})\oplus\cdots\oplus(\mathfrak{o}/\mathfrak{p}^{\lambda_n}),$$
where $\lambda_1\ge\ldots\ge\lambda_n$ are positive integers. 
\end{proposition}

In this case, we say $M$ is a $\mathfrak{o}$-module of type $\lambda=(\lambda_1,\ldots,\lambda_n)$, and we write 
\begin{equation}\label{eq:def_length}
    l(M)=|\lambda|=\lambda_1+\cdots+\lambda_n,
\end{equation}
called the \emph{length} of $M$. It is clear that the isomorphism classes of $\mathfrak{o}$-modules are totally determined by the partition $\lambda$. 

\begin{defi}

Let $\lambda,\mu^{(1)},\ldots,\mu^{(r)}\in\Y$, and let $M$ be a $\mathfrak{o}$-module of type $\lambda$. We define $G_{\mu^{(1)},\ldots,\mu^{(r)}}^\lambda(\mathfrak{o})$ to be the number of chains of submodules of $M$:

$$M=M_0\supset M_1\supset\ldots\supset M_r$$
such that $M_{i-1}/M_i$ has type $\mu^{(i)}$ for $1\le i\le r$.
\end{defi}

\begin{proposition}\label[proposition]{prop: properties of the Hall polynomial}
$G_{\mu^{(1)},\ldots,\mu^{(r)}}^\lambda(\mathfrak{o})$ has the following properties:
\begin{enumerate}
\item $G_{\mu^{(1)},\ldots,\mu^{(r)}}^\lambda(\mathfrak{o})=0$ unless $|\lambda|=|\mu^{(1)}|+\cdots+|\mu^{(r)}|$.

\item $G_{\mu^{(1)},\ldots,\mu^{(r)}}^\lambda(\mathfrak{o})$ does not depend on the order of the subscripts, i.e., for any permutation $\{i_1,\ldots,i_r\}=\{1,2,\ldots,r\}$, we have $G_{\mu^{(1)},\ldots,\mu^{(r)}}^\lambda(\mathfrak{o})=G_{\mu^{(i_1)},\ldots,\mu^{(i_r)}}^\lambda(\mathfrak{o})$.

\item $G_{\mu^{(1)},\ldots,\mu^{(r)}}^\lambda(\mathfrak{o})$ is a ``polynomial in $q$", i.e., there exists a polynomial $g_{\mu^{(1)},\ldots,\mu^{(r)}}^\lambda(t)\in\Z[t]$, independent of $\mathfrak{o}$, such that $G_{\mu^{(1)},\ldots,\mu^{(r)}}^\lambda(\mathfrak{o})=g_{\mu^{(1)},\ldots,\mu^{(r)}}^\lambda(q)$. 
\end{enumerate}
\end{proposition}

\begin{proof}
See (4.3) of \cite[Chapter II]{Macdonald}. 
\end{proof}

With the above preparation, we now state the definition of Hall algebra. 

\begin{defi}\label{def:hall_alg}
We define $H=H(\mf{o})$ to be the free $\Z$-algebra on a basis $(u_\lambda)$ indexed by all partitions $\lambda$ and a product over $H$ given by

$$u_\mu u_\nu=\sum_{\lambda\in\Y}G_{\mu,\nu}^\lambda (\mf{o})u_\lambda,\quad\forall \mu,\nu\in\Y.$$
\end{defi}

$H(\mathfrak{o})$ is a commutative and associative ring with identity element $u_0$. It is generated (as $\Z$-algebra) by the elements $u_{(1^r)}$, and they are algebraically independent over $\Z$. Also, it has a canonical graded ring structure 
$$H(\mathfrak{o})=\bigoplus_{k\ge 0}H^k(\mathfrak{o}),$$
where $H^k(\mathfrak{o})$ is the free additive group with basis $\{u_\lambda\mid\lambda\in\Y, |\lambda|=k\}$.

\subsection{Hecke ring} \label{subsec:hecke_ring}
The Hecke ring structure plays an essential role throughout our method. See \cite[Chapter V]{Macdonald} for further details. We use the term ``Hecke ring'' in this subsection since it is the original terminology from \cite{Macdonald}, though it is also often called the type A spherical Hecke algebra in modern works.

Let $G=\GL_n(F)$ be the group of all invertible $n\times n$ matrices over $F$. Also, let $G^+=G\cap M_{n\times n}(\mathfrak{o})$ be the subsemigroup of $G$ consisting of all matrices $x\in G$ with entries $x_{ij}\in\mathfrak{o}$, and let $K=\GL_{n}(\mathfrak{o})=G^+\cap(G^+)^{-1}$ so that $K$ consists of all $x\in G$ with entries $x_{ij}\in\mathfrak{o}$ and $\det(x)$ a unit in $\mathfrak{o}$.
\begin{defi}\label{defi: Hecke algebra}
Let $L(G,K)$ (resp. $L(G^+,K)$) denote the space of all complex-valued continuous functions of compact support on $G$ (resp. $G^+$) which are bi-invariant with respect to $K$, i.e., such that $f(k_1xk_2)=f(x)$ for all $x\in G$ (resp. $G^+$) and $k_1,k_2\in K$. We define a multiplication on $L(G,K)$ as follows: for $f,g\in L(G,K)$,
$$(f*g)(x)=\int_Gf(xy^{-1})g(y)dy,$$
where $dy$ is the unique Haar measure on $G$ such that $K$ has measure $1$, so that the indicator function on $K$ is the multiplicative identity of $L(G,K)$. This product is associative and commutative, which makes $L(G^+,K)$ a subring of $L(G,K)$. Also, let $H(G,K)$ (resp. $H(G^+,K)$) denote the subspace of $L(G,K)$ (resp. $L(G^+,K)$) consisting of functions with integer values, which shall be called the \emph{Hecke ring} of $G$ (resp. $G^+$). 
\end{defi}
\begin{proposition}
Every $f\in L(G,K)$ (resp. $f\in L(G^+,K)$) could be written as a finite linear combination of the form $c_\mu$, where $\mu=(\mu_1,\ldots,\mu_n),\mu_1\ge\cdots\ge\mu_n\in\Z^n$ (resp. $\mu=(\mu_1,\ldots,\mu_n),\mu_1\ge\cdots\ge\mu_n\in\Z^n_{\ge 0}$), and $c_\mu$ is the characteristic function of the double coset
$$K\diag_{n\times n}(\pi^{\mu_1},\ldots,\pi^{\mu_n})K.$$
The similar statements for $H(G,K),H(G^+,K)$ translate mutatis mutandis.
\end{proposition}

\subsection{Hall-Littlewood polynomials and the ring of symmetric functions}

\begin{defi}\label{defi: symmetric functions}
We write
$$\Lambda_n=\Z[x_1,\ldots,x_n]^{S_n}=\bigoplus_{k\ge 0}\Lambda_n^k$$
for the \emph{ring of symmetric polynomials} over the independent variables $x_1,\ldots,x_n$. Here $\Lambda_n^k$ consists of the homogeneous symmetric polynomials of degree $k$, together with the zero polynomial. For any $m\ge n$, let 
$$\rho_{m,n}:\L_m\rightarrow\L_n$$
be the homomorphism sending every $x_i,1\le i\le n$ to itself, and the other $x_i$ to zero. On restriction to $\L_m^k$ we have homomorphisms
$$\rho_{m,n}^k:\L_m^k\rightarrow\L_n^k$$
for all $k\ge 0$ and $m\ge n$, which are always surjective.
Let $\Lambda^k=\varprojlim\Lambda_n^k$ denote the inverse limit of the $\Z$-modules $\L_n^k$ relative to the homomorphisms $\rho_{m,n}^k$, which is the free $\Z$-module of symmetric functions of degree $k$. Then we call $\Lambda:=\bigoplus_{k\ge 0}\Lambda^k$ the \emph{ring of symmetric functions}. Also, let $\Lambda_\Q:=\Lambda\otimes_\Z\Q,\Lambda[t]:=\Lambda\otimes_\Z\Z[t]$ denote the ring of symmetric function with coefficients in $\Q$ and $\Z[t]$ respectively. Denote the canonical projection maps $\Lambda \to \Lambda_n$ by $\rho_{\infty,n}$.
\end{defi}

\begin{defi}\label{def:Lambda_sq}
Let $\L_n^\sq=\Z[x_1^2,\ldots,x_n^2]^{S_n}$ be the subring of $\Lambda_n$ where the degree of each $x_i$ are required to be even. We define $\Lambda^\sq \subset \Lambda$ to be the subring such that $\rho_{\infty,n}(\Lambda^\sq) = \Lambda_n^\sq$ for every $n$.

% It has a canonical graded ring structure
% $$\L_n^\sq=\bigoplus_{k\ge 0}\L_n^{k,\sq},\quad\L_n^{k,\sq}=\L_n^\sq\cap\L_n^{2k}.$$
% On restriction of $\rho_{m,n}$ in \Cref{defi: symmetric functions} to $\L_m^{k,\sq}$ we have homomorphisms 
% $$\rho_{m,n}^{k,\sq}: \L_m^{k,\sq}\rightarrow\L_n^{k,\sq}$$
% for all $k\ge 0$ and $m\ge n$. Let $\Lambda^{k,\sq}=\varprojlim\Lambda_n^{k,\sq}$ denote the inverse limit of the $\Z$-modules $\L_n^{k,\sq}$ relative to the homomorphisms $\rho_{m,n}^{k,\sq}$, and $\Lambda^{\sq}:=\bigoplus_{k\ge 0}\Lambda^{k,\sq},\Lambda_\Q^\sq:=\Lambda^\sq\otimes_\Z\Q$. 
 
% When we view $\L_n$ as a $\L_n^\sq$-module, we have natural decomposition of graded modules
% $$\L_n=\L_n^\even\oplus\L_n^\odd,\quad\L_n^\even=\bigoplus_{k\ge 0}\L_n^{2k},\L_n^\odd=\bigoplus_{k\ge 0}\L_n^{2k+1}.$$
\end{defi}

One should view $\Lambda$ (resp. $\L^\sq$) heuristically as a ring of symmetric polynomials in countably many independent variables $x_1,x_2,\ldots$ (resp. $x_1^2,x_2^2,\ldots$). There are many natural bases for the ring $\L$ of symmetric functions, one being the \emph{Hall-Littlewood symmetric functions} defined as follows.

\begin{defi}
Let $\lambda=(\l_1,\ldots,\l_n)\in\Y$ be a partition of length $\le n$. Then, the \emph{Hall-Littlewood polynomial} $P_\lambda(x_1,\ldots,x_n;t)$ is defined by

\begin{equation}\label{eq: Hall Littlewood polynomial}
P_\l(x_1,\ldots,x_n;t) = \frac{1}{V_\l(t)} \sum_{\sigma \in S_n} \sigma\left(x_1^{\lambda_1}\cdots x_n^{\lambda_n} \prod_{1 \leq i < j \leq n} \frac{x_i-tx_j}{x_i-x_j}\right),
\end{equation}
where $V_\l(t) = \prod_{i \in \Z_{\ge 1}} \frac{(1-t)\cdots(1-t^{m_i(\lambda)})}{(1-t)^{m_i(\l)}}=\frac{1}{(1-t)^n}\prod_{i\in\Z_{\ge 1}}(t;t)_{m_i(\lambda)}$. Here the notation $(a;q)_n := (1-a)(1-aq) \cdots (1-aq^{n-1}),n\ge 0$ is the $q$-Pochhammer symbol, with $(a;q)_0 =1$.
\end{defi}

\begin{proposition}
The Hall-Littlewood polynomials $P_\l(x_1,\ldots,x_n;t)$ satisfy the following properties:

\begin{enumerate}
\item They have the form
$$P_\l(x_1,\ldots,x_n;t) = x_1^{\l_1}x_2^{\l_2}\cdots x_n^{\l_n} + \text{(lower-order monomials in the lexicographic order)}.$$

\item When $\lambda\in\Y$ ranges over all partitions of length $\le n$, the Hall-Littlewood polynomials $P_\l(x_1,\ldots,x_n;t)$ form a $\Z[t]$-basis of $\L_n[t]$, where $\L_n[t]:=\Z[t][x_1,\ldots,x_n]^{S_n}$ is the ring of symmetric polynomials in $n$ variables $x_1,\ldots,x_n$ with coefficients in $\Z[t]$.

\item\label{item: finite to infinite} Let $\lambda=(\l_1,\ldots,\l_n)\in\Y$ be a partition of length $\le n$. Then $P_\l(x_1,\ldots,x_n;t)=P_\l(x_1,\ldots,x_n,0;t)$.
\end{enumerate}
\end{proposition}

By \eqref{item: finite to infinite} we may pass to the limit and define $P_\lambda(x_1,x_2,\ldots;t)$ to be the element of $\Lambda[t]$ whose image in $\Lambda_n[t]$ for each $n\ge l(\lambda)$ is $P_\lambda(x_1,\ldots,x_n;t)$. We call $P_\lambda(x_1,x_2,\ldots;t)$ the \emph{Hall-Littlewood symmetric function} corresponding to the partition $\lambda$, which is homogeneous of degree $|\lambda|$. We define $P_\l(x_1,x_1t,x_2,x_2t,\ldots;t)$ in the same way, as the elements of $\La$ which map to $P_\l(x_1,x_1t,\ldots,x_n,x_nt;t)$ under the natural map to $\La_n$, and similarly for $P_\l(x_1^2,x_2^2,\ldots;t^2)$. 

The dual Hall-Littlewood polynomials $Q_\l(x_1,x_2,\ldots;t)$, which form a dual basis to the $P_\l(x_1,x_2,\ldots;t)$ under a natural inner product defined in \cite[Chapter III.4]{Macdonald}, are given by 
\begin{equation}\label{eq:def_Q}
    Q_\l(x_1,x_2,\ldots;t) = \prod_{i\ge 1}(t;t)_{m_i(\lambda)}P_\l(x_1,x_2,\ldots;t).
\end{equation}
We sometimes omit the parameter $t$ when it is clear from the context. Because the $P_\l,\lambda\in\Y$ form a basis for the vector space of symmetric polynomials in $n$ variables, there exist symmetric polynomials $P_{\l/\mu}(x_1,\ldots,x_k;t)\in\Lambda_k[t]$ indexed by $\l\in\Y_{n+k},\mu\in\Y_n$ which are defined by
\begin{equation*}\label{eq: skew Hall Littlewood polynomial}
P_\l(x_1,\ldots,x_{n+k};t) = \sum_{\mu \in \Y_n} P_{\l/\mu}(x_1,\ldots,x_k;t) P_\mu(x_{k+1},\ldots,x_{n+k};t).
\end{equation*}
Similarly, for $\l\in\Y_{n+k},\mu\in\Y_n$ and $k \geq 1$ arbitrary, define $Q_{\l/\mu}(x_1,\ldots,x_k;t)$ by
\begin{equation*}
Q_\l(x_1,\ldots,x_{n+k};t) = \sum_{\mu \in \Y_n} Q_{\l/\mu}(x_1,\ldots,x_k;t) Q_\mu(x_{k+1},\ldots,x_{n+k};t).
\end{equation*}
The polynomials $P_{\l/\mu},Q_{\l/\mu}$ are called \emph{skew Hall-Littlewood polynomials}. For any $\l,\mu\in\Y$ and large enough $k$, we have
$$P_{\l/\mu}(x_1,\ldots,x_k,0;t)=P_{\l/\mu}(x_1,\ldots,x_k;t),$$
$$Q_{\l/\mu}(x_1,\ldots,x_k,0;t)=Q_{\l/\mu}(x_1,\ldots,x_k;t).$$
Hence there are symmetric functions $P_{\l/\mu}(x_1,x_2,\ldots;t),Q_{\l/\mu}(x_1,x_2,\ldots;t)\in\L[t]$ which map to $P_{\l/\mu}(x_1,\ldots,x_k;t),Q_{\l/\mu}(x_1,\ldots,x_k;t)$ under the maps $\L\mapsto\L_k$ respectively. In particular, we have $Q_{\l/\mu}(x_1,x_2,\ldots;t)=\prod_{i\ge 1}\frac{(t;t)_{m_i(\l)}}{(t;t)_{m_i(\mu)}}P_{\l/\mu}(x_1,x_2,\ldots;t)$, and
\begin{equation}\label{eq:Qs_agree}
    Q_{\l/(0)}(x_1,x_2,\ldots;t) = Q_\l(x_1,x_2,\ldots;t).
\end{equation}

\begin{lemma}\label[lemma]{lem:asym_cauchy}
(Skew Cauchy identity) Let $\mu,\nu\in \Y$, and $x_1,x_2,\ldots, y_1,y_2,\ldots$ be indeterminates. Then
\begin{multline}\label{eq:asym_cauchy}
    \sum_{\kappa \in \Y} Q_{\kappa/\mu}(y_1,y_2,\ldots;t) P_{\kappa/\nu}(x_1,x_2,\ldots;t) \\
    = \Pi_t(x_1,x_2,\ldots; y_1,y_2,\ldots)\sum_{\l\in \Y} Q_{\nu/\l}(y_1,y_2,\ldots;t) P_{\mu/\l}(x_1,x_2,\ldots;t), 
\end{multline}
where the Cauchy kernel $\Pi_t(x_1,x_2,\ldots; y_1,y_2,\ldots)$ has the form 
\begin{equation}\label{eq: Cauchy kernel}
    \Pi_t(x_1,x_2,\ldots; y_1,y_2,\ldots) = \prod_{i,j} \frac{1-tx_iy_j}{1-x_iy_j}=\exp\left(\sum_{n\ge1}\frac{1-t^n}{n}\sum_{i\ge 1}x_i^n\cdot\sum_{j\ge 1}y_j^n\right),
\end{equation}
and \eqref{eq:asym_cauchy} is interpreted as an equality of formal power series in the variables.
\end{lemma}

When we take $t$ to be a real number, for any sequence of real numbers $\rho=(a_1,a_2,\ldots)$, we define
\begin{equation}
    P_{\la/\mu}(\rho;t) = \lim_{n \to \infty} P_{\la/\mu}(a_1,\ldots,a_n;t)
\end{equation}
if the limit exists.

% \begin{defi}
% Fix the parameter $t$ to be a real number. For any sequence of real numbers $\rho=(a_1,a_2,\ldots)$, we define
% \begin{equation}
%     P_{\la/\mu}(\rho;t) = \lim_{n \to \infty} P_{\la/\mu}(a_1,\ldots,a_n;t)
% \end{equation}
% if the limit exists. We say that the sequence $\rho$ is \emph{Hall-Littlewood nonnegative} or simply \emph{nonnegative} if it takes non-negative values on the skew Hall-Littlewood symmetric functions: $P_{\lambda/\mu}(\rho;t)\ge 0$ for all $\lambda,\mu\in\Y$.
% \end{defi}

% \begin{rmk}
%     When $t \in [0,1)$, it follows from \cite[Theorem 1.4]{matveev2019macdonald} that the sequence $\rho$ is Hall-Littlewood nonnegative if and only if $a_i \geq 0$ for all $i$ and $\sum_{i \geq 1} a_i < \infty$. Nevertheless, for negative $t$ (which we must use in the Hermitian case), sequences such as $1,t,t^2,\ldots$ will contain negative elements but still be Hall-Littlewood nonnegative.
% \end{rmk}

\begin{defi}
Fix $t$ to be real. Let $\rho_1=(a_1,a_2,\ldots),\rho_2=(b_1,b_2,\ldots)$ be sequences of real numbers such that
$$\Pi_t(\rho_1;\rho_2)=\sum_\l P_\l(\rho_1;t)Q_\l(\rho_2;t)=\exp\left(\sum_{n\ge1}\frac{1-t^n}{n}p_n(\rho_1)p_n(\rho_2)\right)<\infty,$$
where $p_n$ is the $n$th power sum, and each summand $P_\l(\rho_1;t)Q_\l(\rho_2;t)$ is nonnegative. Then the corresponding \emph{Hall-Littlewood measure} is the probability measure on $\Y$ defined by
$$\mathbf{P}(\lambda)=\frac{P_\l(\rho_1;t)Q_\l(\rho_2;t)}{\Pi_t(\rho_1;\rho_2)},\quad\forall\l\in\Y.$$
\end{defi}

When all the $a_i$ are in geometric progression with common ratio $t$, the skew Hall-Littlewood polynomials take a simple form due to \cite{kirillov1998new}:

\begin{theorem}
(Principal specialization formula) For $\lambda,\mu\in\Y$, we have
\begin{equation}\label{eq: principal specialization formula for P}
P_{\l/\mu}(u,ut,\ldots;t)=u^{|\l|-|\mu|}t^{n(\l/\mu)}\prod_{i\ge 1}\frac{(t^{1+\l'_i-\mu_i'};t)_{m_i(\mu)}}{(t;t)_{m_i(\l)}},
\end{equation}
\begin{equation}\label{eq: principal specialization formula for Q}
Q_{\l/\mu}(u,ut,\ldots;t)=u^{|\l|-|\mu|}t^{n(\l/\mu)}\prod_{i\ge 1}\frac{(t^{1+\l'_i-\mu_i'};t)_{m_i(\mu)}}{(t;t)_{m_i(\mu)}},
\end{equation}
where $n(\l/\mu)=\sum_{i\ge 1}\begin{pmatrix}\l_i'-\mu_i' \\ 2\end{pmatrix}$. In particular, when $\mu$ is the zero partition, we have
\begin{equation}\label{eq: non skew principal specialization formula for P}
P_\l(u,ut,\ldots;t)=u^{|\l|}t^{n(\l)}\prod_{i\ge 1}\frac{1}{(t;t)_{m_i(\l)}},
\end{equation}
\begin{equation}\label{eq: non skew principal specialization formula for Q}
Q_\l(u,ut,\ldots;t)=u^{|\l|}t^{n(\l)},
\end{equation}
where $n(\l)=\sum_{i\ge 1}\begin{pmatrix}\l_i' \\ 2\end{pmatrix} = \sum_{i \geq 1} (i-1) \la_i$.
\end{theorem}

The below theorem connects the story of Hall algebra, Hecke ring, and the ring of symmetric polynomials.

\begin{theorem}\label{thm: left half}
Let $G^+=\GL_n(F)\cap M_n(\mathfrak{o}),K=\GL_n(\mathfrak{o})$. We have canonical homomorphisms
\begin{equation}
H(\mathfrak{o})\overset{(1)}{\longrightarrow}H(G^+,K)\overset{(2)}{\longrightarrow}\Lambda_n\otimes_\Z\Q,
\end{equation}
where
\begin{itemize}
    \item $(1)$ is surjective and given by $u_\lambda\mapsto c_\lambda$, and the kernel is generated by the $u_\lambda$ such that $l(\lambda)>n$.
    \item $(2)$ is an isomorphism given by $c_\lambda\mapsto q^{-n(\lambda)}P_\lambda(x_1,\ldots,x_n;q^{-1})$.
\end{itemize} 
\end{theorem}

\begin{proof}
The map $(1)$ is given in (2.6) of \cite[Chapter V]{Macdonald}, and the map $(2)$ is given in (2.7) of \cite[Chapter V]{Macdonald}.
\end{proof}

Taking the inverse limit (viewed as graded rings) of the maps given in \Cref{thm: left half} brings us to the following elegant and classical description of the Hall algebra.

\begin{corollary}\label[corollary]{cor: Hall and symmetric functions}
Let $\psi:H(\mathfrak{o})\otimes_\Z\Q\rightarrow\L_\Q$ be the $\Q$-linear mapping defined by
\begin{equation}\label{eq: standard form of isomorphism}
\psi(u_\l)=q^{-n(\l)}P_\l(x_1,x_2,\ldots;q^{-1}).
\end{equation}
Then $\psi$ is an isomorphism of rings.
\end{corollary}

Our results in the alternating and Hermitian cases, given in \Cref{thm: Hall and Symmetric functions}, are of the same ilk and should be viewed in the context of this one.

\section{The Alternating Hall Module} \label{sec:alternating}

\subsection{Alternating module background}

There is a natural homomorphism from the Hall algebra $H=H(\mathfrak{o})$ to the Hecke ring $H(G^+,K)$, as shown in \Cref{thm: left half}. As an analogy, one expects to define an ``alternating Hall module'' $H^{\alt}=H^{\alt}(\mathfrak{o})$ over the Hall algebra $H(\mathfrak{o})$ and connect it with the alternating Hecke module $H(G^{+\alt},K)$, which is the goal of this section.

Let $M$ be a finite $\mathfrak{o}$-module, equipped with a bilinear alternating pairing $\langle\cdot,\cdot\rangle:M\times M\rightarrow F/\mathfrak{o}$. Here recall that alternating means $\langle x,y \rangle = - \langle y, x \rangle$. For any $\mathfrak{o}$-submodule $N$ of $M$, we write $\langle\cdot,\cdot\rangle\mid_N$ to denote the restriction of the alternating pairing to $N$. Let $\hat M=\text{Hom}(M,F/\mathfrak{o})$ equipped with the canonical $\mathfrak{o}$-module structure. The following properties are equivalent for the alternating pairing:
\begin{enumerate}
\item The map $M\rightarrow\hat M, x\mapsto\langle x,\cdot\rangle$ gives an isomorphism from $M$ to $\hat M$. \label{item: alt_iso}

\item If $x\in M$ satisfies $\langle x,y\rangle=0$ for every $y\in M$, then $x=0$.\label{item: alt_x=0}
\end{enumerate}
\begin{defi}
If the above properties \eqref{item: alt_iso} and \eqref{item: alt_x=0}
are satisfied, we say the pairing is \emph{regular}, and $M$ has an \emph{alternating $\mathfrak{o}$-module} structure. 
\end{defi}

Let $(M,\langle\;,\;\rangle),(M',\langle\;,\;\rangle')$ be two alternating $\mathfrak{o}$-modules. We say these two $\mathfrak{o}$-modules are \emph{isomorphic}, if there exists an $\mathfrak{o}$-module isomorphism $\varphi:M\rightarrow M'$ that preserves the pairing, i.e., for all $x,y\in M$,
$$\langle x,y\rangle=\langle\varphi(x),\varphi(y)\rangle'.$$
\begin{theorem}
Every alternating $\mathfrak{o}$-module $M$ is isomorphic to one of the form
\begin{equation}
M\cong(\mathfrak{o}/\mathfrak{p}^{\lambda_1}\oplus\mathfrak{o}/\mathfrak{p}^{\lambda_1})\oplus\cdots\oplus(\mathfrak{o}/\mathfrak{p}^{\lambda_n}\oplus\mathfrak{o}/\mathfrak{p}^{\lambda_n})
\end{equation}
for some $\lambda_1\ge\ldots\ge\lambda_n\ge 1$ with the pairing $\langle\cdot,\cdot\rangle$ given by $\langle x,y\rangle= x^T(\pi_\lambda^{\alt})^{-1}y$ for any $x=(x_1,\ldots,x_{2n}),y=(y_1,\ldots,y_{2n})\in M$, where
\begin{equation}
\pi^{\alt}_\lambda=\diag\left({\begin{pmatrix} 0 & \pi^{\lambda_1} \\ -\pi^{\lambda_1} & 0\end{pmatrix}},\ldots,{\begin{pmatrix} 0 & \pi^{\lambda_n} \\ -\pi^{\lambda_n} & 0\end{pmatrix}}\right).
\end{equation}
\end{theorem}
\begin{proof}
See the corollary of \cite[Theorem 4.1]{Symplectic_modules}.
\end{proof}

In this case, we say $M$ is an alternating $\mathfrak{o}$-module of type $\lambda=(\lambda_1,\ldots,\lambda_n)$. The isomorphism class of alternating $\mathfrak{o}$-modules is totally determined by the partition $\lambda$.

For any $\mathfrak{o}$-submodule $N$ of the alternating $\mathfrak{o}$-module $M$, let
\begin{equation}\label{eqref: alt_orthogonal submodule}
N^\perp=\{x\in M\mid\langle x,y\rangle=0,\quad\forall y\in N\}\end{equation}
be the \emph{orthogonal submodule of $N$ in $M$}. In this case, the $\mathfrak{o}$-module homomorphism
$$M\rightarrow\hat N:=\Hom(N,F/\mathfrak{o}),\quad x\mapsto\langle x,\cdot\rangle$$
is surjective with kernel $N^\perp$. Note that the isomorphism $N\cong\hat N$ always holds. This can be seen by observing that the mapping commutes with direct sums, and therefore, it is sufficient to verify the case when $N$ is cyclic, which is straightforward. Therefore, we obtain the following:

\begin{proposition}\label[proposition]{thm:perp_alt_props}
If $M$ is an alternating $\mathfrak{o}$-module and $N,N'\subset M$, then
\begin{enumerate}
\item $N\cong M/N^\perp,N^\perp\cong M/N$ as $\mathfrak{o}$-modules.

\item $\#N\cdot\#N^\perp=\#M$.

\item $N^{\perp\perp}=N$.

\item If $N$ is regular (i.e., the restriction of the pairing over $N$ is regular), then $N^\perp$ is also regular, and $M\cong N\oplus N^\perp$.

\item $N\subset N'\Leftrightarrow N'^\perp\subset N^\perp$. \label{item: alt_include perp}
\end{enumerate}
\end{proposition}

\begin{proof}
See \cite[Proposition 2.2]{Symplectic_modules}.
\end{proof}

An $\mathfrak{o}$-submodule $N\subset M$ is \emph{isotropic} if $N\subset N^\perp$, i.e. $\langle N,N\rangle=0$. In this case, the alternating pairing over $M$ induces a regular alternating pairing over $N^\perp/N$, which gives this $\mathfrak{o}$-module an alternating structure. Also, by part \eqref{item: alt_include perp} of the above theorem, every $\mathfrak{o}$-submodule $P\subset N$ is also isotropic.

\subsection{The alternating Hecke module $H(G^{\alt},K)$}

The goal of this subsection is to provide a  background for the alternating Hecke module studied in \cite{shen2024non}, which turns out in \Cref{thm: alt_Hall and Hecke} to share the same structure coefficient as the alternating Hall module $H^{\alt}(\mathfrak{o})$.

Let $G=\GL_{2n}(F)$ be the group of all invertible $2n\times 2n$ matrices over $F$. Also, let $G^+=G\cap M_{2n\times 2n}(\mathfrak{o})$ be the subsemigroup of $G$ consisting of all matrices $x\in G$ with entries $x_{ij}\in\mathfrak{o}$, and let $K=\GL_{2n}(\mathfrak{o})=G^+\cap(G^+)^{-1}$ so that $K$ consists of all $x\in G$ with entries $x_{ij}\in\mathfrak{o}$ and $\det(x)$ a unit in $\mathfrak{o}$.

Let $L(G,K)$ (resp. $L(G^+,K)$) denote the space of all complex-valued continuous functions of compact support on $G$ (resp. $G^+$) which are bi-invariant with respect to $K$, i.e., such that $f(k_1xk_2)=f(x)$ for all $x\in G$ (resp. $G^+$) and $k_1,k_2\in K$. Let $H(G,K)$ (resp. $H(G^+,K)$) denote the subspace of $L(G,K)$ (resp. $L(G^+,K)$) consisting of functions with integer values. The convolution over these spaces provides a ring structure, as previously defined in \Cref{defi: Hecke algebra}. 

\begin{defi}
For all $x\in G$, we write $x^T$ as the transpose of $x$, and $x^{-T}=(x^{-1})^{T}$ as the inverse transpose of $x$. Let $G^{\alt}=G\cap\Alt_{2n}(F)=\{x\in G\mid x^T=-x\}$ denote the set of invertible alternating matrices, and $G^{+{\alt}}=G^+\cap\Alt_{2n}(\mathfrak{o})$ denote the subset of $G$ consisting of matrices with entries in $\mathfrak{o}$.
\end{defi}

\begin{defi}
Let $L(G^{\alt},K)$ (resp. $L(G^{+{\alt}},K)$) denote the space of all complex-valued continuous functions of compact support on $G^{\alt}$ (resp. $G^{+\alt}$) which are invariant with respect to $K$, i.e., such that
$$f(k^Txk)=f(x)$$
for all $x\in G$ (resp. $x\in G^+$) and $k\in K$. We may regard $L(G^{+{\alt}},K)$ as a subset of $L(G^{\alt},K)$. 
\end{defi}

We define a multiplication of $L(G,K)$ over $L(G^{\alt},K)$, giving the latter an $L(G,K)$-module structure, as follows. For $g\in L(G,K)$, $f\in L(G^{\alt},K)$,
$$(g*f)(x):=\int_Gg(y)f(y^{-1}xy^{-T})dy.$$
Since $g$ and $f$ are compactly supported, the integration is over a compact set. This product satisfies $g_1*(g_2*f)=(g_1*g_2)*f$ and $c_0*f=f$, hence defines a $L(G,K)$ (resp. $L(G^+,K)$)-module structure of $L(G^{\alt},K)$ (resp. $L(G^{+\alt},K)$). We use the same notation $*$ for the convolution of the Hecke ring itself and the convolution of the Hecke ring over the alternating Hecke module, but it is easy to distinguish these two based on context.

Each function $f\in L(G^{\alt},K)$ is constant on the orbits $\{k^Txk|k\in K\}$ in $G^{\alt}$. These orbits are compact, open, and mutually disjoint. Since $f$ has compact support, it follows that $f$ takes non-zero values on only finitely many orbits $\{k^Txk|k\in K\}$, and hence can be written as a finite linear combination of their characteristic functions. Therefore, the characteristic functions of these orbits in $G^{\alt}$ form a $\C$-basis of $L(G^{\alt},K)$.

\begin{defi}\label{defi: alt_integer value Hecke module}
Let $H(G^{\alt},K)$ (resp. $H(G^{+{\alt}},K)$) denote the subspace of $L(G^{\alt},K)$ (resp. $L(G^{+{\alt}},K)$) consisting of all integer-valued functions, which has a $H(G,K)$ (resp. $H(G^+,K)$)-module structure as above. We may regard $H(G^{+{\alt}},K)$ as a subset of $H(G^{\alt},K)$.
\end{defi}

Clearly we have
$$L(G^{\alt},K)=H(G^{\alt},K)\otimes_\Z\C,L(G^{+{\alt}},K)=H(G^{+{\alt}},K)\otimes_\Z\C.$$
Consider an orbit $\{k^Txk|k\in K\}$, where $x\in G^{\alt}$. By multiplying $x$ by a suitable power of $\pi$ (the generator of $\mathfrak{p}$) we can bring $x$ into $G^{+{\alt}}$. As stated on page 483 of \cite{Hironaka}, each orbit $\{k^Txk|k\in K\}$ has a unique representative of the form
$$\pi^{\alt}_\lambda=\text{diag}\left({\begin{pmatrix} 0 & \pi^{\lambda_1} \\ -\pi^{\lambda_1} & 0\end{pmatrix}},\ldots,{\begin{pmatrix} 0 & \pi^{\lambda_n} \\ -\pi^{\lambda_n} & 0\end{pmatrix}}\right),$$
where $\lambda_1\ge\lambda_2\ge\ldots\ge\lambda_n$. We have $\lambda_n\ge 0$ if and only if $x\in G^{+{\alt}}$. 

Let $c^{\alt}_\lambda$ denote the characteristic function of the orbit $\{k^T\pi^{\alt}_\lambda k|k\in K\}$. Then we have the $c^{\alt}_\lambda$ (resp. the $c^{\alt}_\lambda$ such that $\lambda_n\ge 0$) form a $\Z$-basis of $H(G^{\alt},K)$(resp. $H(G^{+{\alt}},K)$). The characteristic function $c_0$ of $K$ is the identity element of $H(G,K)$ and $H(G^+,K)$, and this also plays the role of the identity element when multiplying with elements in $H(G^{\alt},K)$ and $H(G^{+{\alt}},K)$. 

Let $\mu\in\Y_{2n},\nu\in\Y_n$ be partitions. Then, the product $c_\mu *c^{\alt}_\nu$ has the form
\begin{equation}\label{eq:structure coefficient_alt}c_\mu *c^{\alt}_\nu=\sum_{\la \in \Y_n} g_{\mu,\nu}^{\alt,\lambda}(q)c^{\alt}_\lambda,
\end{equation}
where \cite[Corollary 3.4]{shen2024non} proves that $g_{\mu,\nu}^{\alt,\lambda}(q)\in\Z[q]$ is a polynomial in $q$ independent of $\mathfrak{o}$. Explicitly, we have
$$g_{\mu,\nu}^{\alt,\lambda}(q)=(c_\mu*c_\nu^{\alt})(\pi_\lambda^{\alt})=\int_G c_\mu(y)c_\nu^{\alt}(y^{-1}\pi_\lambda^{\alt} y^{-T})dy.$$
Since $c_\mu(y)$ vanishes for $y$ outside $K\pi_\mu K$, the integration is over this orbit, which we shall write as a disjoint union of left cosets, say
\begin{equation}\label{eq: alt_disjoint left coset}K\pi_\mu K=\bigsqcup_j y_jK\quad (y_j\in K\pi_\mu ).\end{equation}
This gives 
\begin{equation}\label{eq: alt_sum of left cosets}g_{\mu,\nu}^{\alt,\lambda}(q)=\sum_j\int_{y_jK}c^{\alt}_\nu(y^{-1}\pi_\lambda^{\alt} y^{-T})dy=\sum_jc^{\alt}_\nu(y_j^{-1}\pi_\lambda^{\alt} y_j^{-T}).\end{equation}

\subsection{The alternating Hall polynomial}

The following theorem is the critical step that connects the alternating Hall module and the alternating Hecke module. 

\begin{theorem}\label{thm: alt_Hall and Hecke}
Let $\mu\in\Y_{2n}$ and let $\lambda,\nu\in\Y_n$. Let $M$ be an alternating module of type $\lambda$. Then the polynomial $g_{\mu,\nu}^{\alt,\lambda}(q)$ is equal to $G_{\mu,\nu}^{\alt,\lambda}(\mathfrak{o})$, the number of submodules $M'\subset M$ such that $M/M'$ is an $\mathfrak{o}$-module of type $\mu$, $M'^\perp$ is isotropic, and $M'/M'^\perp$ is an alternating $\mathfrak{o}$-module of type $\nu$.
\end{theorem}

\begin{proof} 
Consider the lattice $L=\bigoplus_{i=1}^{2n}\mathfrak{o}$ equipped with an alternating pairing $L\times L\rightarrow F/\mathfrak{o}: \langle x,y\rangle= x^T(\pi_\lambda^{\alt})^{-1}y$, where $x,y\in L$. We view $M$ as the quotient $M=L/\pi_\lambda^{\alt}L$ with bilinear pairing inherited from $L$. It is clear that the left cosets $y_jK$ in \eqref{eq: alt_disjoint left coset} are in one-to-one correspondence with the submodules $M'\subset M$ such that $M/M'$ is a finite $\mathfrak{o}$-module with type $\mu$, under the map $y_j K \mapsto M'=y_jL/\pi_\lambda^{\alt}L$. Then we claim the following are equivalent:
\begin{enumerate}[left=0pt]
\item $c_\nu(y_j^{-1}\pi_\lambda^{\alt} y_j^{-T})=1$.\label{item: alt_c=1}

\item There exists $k_j\in K$ such that $y_j^{-1}\pi_\lambda^{\alt} y_j^{-T}=k_j\pi_\nu^{\alt}k_j^T$.\label{item: alt_existence}

\item $M'\supset M'^\perp$, and $M'/M'^\perp$ is an alternating module of type $\nu$. Here $M'^\perp$ is the orthogonal submodule of $M'$ defined in \eqref{eqref: alt_orthogonal submodule}. Note that it follows from the previous sentence that $M'^\perp$ is isotropic.
\label{item: alt_type}
\end{enumerate}

Now let us prove the above equivalence:

\eqref{item: alt_c=1}$\Rightarrow$\eqref{item: alt_existence}: Trivial.

\eqref{item: alt_existence}$\Rightarrow$\eqref{item: alt_type}: In this case, we have $M'=y_jL/\pi_\lambda^{\alt}L$. We claim that $M'^\perp=y_jk_j\pi_\nu^{\alt}k_j^TL/\pi_\lambda^{\alt}L.$ An element $y_j a \in M'$ will have 
\begin{equation*}
    \langle \pi_\la^{\alt} y_j^{-T} b, y_j a \rangle = b^T y_j^{-1} \pi_\la^{\alt} (\pi_\la^{\alt})^{-1} y_j a = b^T a,
\end{equation*}
which is $0$ in $F/\mf{o}$ for all $a,b\in\bigoplus_{i=1}^{2n}\mathfrak{o}$, hence $\pi_\lambda^{\alt} y_j^{-T}L/\pi_\lambda^{\alt} L \subset M'^\perp$.

Conversely, if $b\in\bigoplus_{i=1}^{2n} F$ satisfies $b^Ta\in\mathfrak{o}$ for all $a\in \bigoplus_{i=1}^{2n}\mathfrak{o}$, then $b\in\bigoplus_{i=1}^{2n}\mathfrak{o}$, which shows the opposite inclusion. Thus
\begin{equation}\label{eq:mperp_alt}
M'^\perp=\pi_\lambda^{\alt} y_j^{-T}L/\pi_\lambda^{\alt}L=y_jk_j\pi_\nu^{\alt}k_j^TL/\pi_\lambda^{\alt}L.
\end{equation} 
where the second equality comes from the hypothesis \eqref{item: alt_existence}. The right hand side of \eqref{eq:mperp_alt} is contained in $M'=y_jL/\pi_\lambda^{\alt}L$, so $M'\supset M'^\perp$. Finally, the alternating pairing restricted to $M'$ is given by \begin{equation}\langle y_jk_jx,y_jk_jy\rangle=x^T(\pi_\nu^{\alt})^{-1}y
\end{equation} 
by \eqref{item: alt_existence}, which shows that $M'/M'^\perp$ is an alternating module of type $\nu$.

\eqref{item: alt_type}$\Rightarrow$\eqref{item: alt_c=1}: In this case, the type of the alternating module $M'/M'^\perp$ implies the alternating pairing restricted over $M'$ has type $\nu$, i.e., $c_\nu(y_j^{-1}\pi_\lambda^{\alt} y_j^{-T})=1$. 

Therefore, following the expression in \eqref{eq: alt_sum of left cosets}, the number of $y_j$ that satisfies the above equivalent properties \eqref{item: alt_c=1}, \eqref{item: alt_existence}, and \eqref{item: alt_type} is $G_{\mu,\nu}^{\alt,\lambda}(\mathfrak{o})$. This ends the proof.
\end{proof}

\Cref{thm: alt_Hall and Hecke} motivates us to generalize the coefficient $G_{\mu,\nu}^{\alt,\lambda}(\mathfrak{o})=g_{\mu,\nu}^{\alt,\lambda}(q)$ which appeared in \eqref{eq:structure coefficient_alt}, based on the alternating Hall module we defined. 

\begin{defi}

Let $\mu^{(1)},\ldots,\mu^{(r)}\in\Y_{2n},\lambda,\nu\in\Y_n$. Let $M$ be an alternating module of type $\lambda$. We define $G_{\mu^{(1)},\ldots,\mu^{(r)},\nu}^{\alt,\lambda}(\mathfrak{o})$ to be the number of chains of submodules of $M$:
$$M=M_0\supset M_1\supset\ldots\supset M_r$$
such that $M_{i-1}/M_i$ has type $\mu^{(i)}$ as an $\mathfrak{o}$-module, $M_i^\perp$ (meaning the orthogonal submodule of $M_i$ inside $M$) is isotropic for $1\le i\le r$, and $M_r/M_r^\perp$ is of type $\nu$. %\RC{the notation $N^\perp$ implicitly depends on the module $N$ is a submodule of; I think here you're considering $M_i^\perp$ with respect to the inclusion $M_i \subset M_{i-1}$ but you should make this clear}  

\end{defi}
\begin{proposition}
$G_{\mu^{(1)},\ldots,\mu^{(r)},\nu}^{\alt,\lambda}(\mathfrak{o})$ has the following properties:
\begin{enumerate}
\item $G_{\mu^{(1)},\ldots,\mu^{(r)},\nu}^{\alt,\lambda}(\mathfrak{o})=0$ unless $|\lambda|=|\mu^{(1)}|+\cdots+|\mu^{(r)}|+|\nu|$.

\item $G_{\mu^{(1)},\ldots,\mu^{(r)},\nu}^{\alt,\lambda}(\mathfrak{o})=\sum_\mu G_{\mu^{(1)},\ldots,\mu^{(r)}}^{\mu}(\mathfrak{o})G_{\mu,\nu}^{\alt,\lambda}(\mathfrak{o})$, where the sum covers $\mu\in\Y$ with $l(\mu)\le 2n$.\label{item: alt and usual polynomial} 

\item $G_{\mu^{(1)},\ldots,\mu^{(r)},\nu}^{\alt,\lambda}(\mathfrak{o})$ does not depend on the order of $\mu^{(1)},\ldots,\mu^{(r)}$, i.e., for any permutation $\{i_1,\ldots,i_r\}=\{1,2,\ldots,r\}$, we have $G_{\mu^{(1)},\ldots,\mu^{(r)},\nu}^{\alt,\lambda}(\mathfrak{o})=G_{\mu^{(i_1)},\ldots,\mu^{(i_r)},\nu}^{\alt,\lambda}(\mathfrak{o})$.

\item $G_{\mu^{(1)},\ldots,\mu^{(r)},\nu}^{\alt,\lambda}(\mathfrak{o})$ is a ``polynomial in $q$", i.e., there exists a polynomial $g_{\mu^{(1)},\ldots,\mu^{(r)},\nu}^{\alt,\lambda}(t)\in\Z[t]$, independent of $\mathfrak{o}$, such that $G_{\mu^{(1)},\ldots,\mu^{(r)},\nu}^{\alt,\lambda}(\mathfrak{o})=g_{\mu^{(1)},\ldots,\mu^{(r)},\nu}^{\alt,\lambda}(q)$. 
\end{enumerate}
\end{proposition}

\begin{proof}
\begin{enumerate}
\item This is because for all $M=M_0\supset M_1\supset\ldots\supset M_r$ such that $M_{i-1}/M_i$ has type $\mu^{(i)}$ as an $\mathfrak{o}$-module, $M_i^\perp$ is isotropic for $1\le i\le r$, and $M_r/M_r^\perp$ is of type $\nu$ as an alternating $\mathfrak{o}$-module, we have 
\begin{align*}
2|\lambda|=l(M_0)=\sum_{i=1}^{r}(l(M_{i-1}/M_i)+l(M_i^\perp/M_{i-1}^\perp))+l(M_r/M_r^\perp)=2|\mu^{(1)}|+\cdots+2|\mu^{(r)}|+2|\nu|,
\end{align*}
where $l(\cdot)$ is the length defined in \eqref{eq:def_length}. 

\item To count the number of chains of modules $M=M_0\supset M_1\supset\ldots\supset M_r$ that satisfy our requirement, we may start from finding a submodule $M_r\subset M$ such that $M_r^\perp$ is isotropic, and $M_r/M_r^\perp$ is an alternating $\mathfrak{o}$-module of type $\nu$. Then we choose the remaining elements of the chain $M=M_0\supset M_1\supset\ldots\supset M_r$ such that $M_{i-1}/M_i$ has type $\mu^{(i)}$ as an $\mathfrak{o}$-module. Since $M_r^\perp$ is isotropic and $M_i^\perp \subset M_r^\perp$, $M_i^\perp$ is automatically isotropic for all $1\le i\le r-1$ by \Cref{thm:perp_alt_props}. Summing over all the types of $M/M_r$ as an $\mathfrak{o}$-module brings the result.

\item We have
\begin{align}\begin{split}
G_{\mu^{(1)},\ldots,\mu^{(r)},\nu}^{\alt,\lambda}(\mathfrak{o})&=\sum_\mu G_{\mu^{(1)},\ldots,\mu^{(r)}}^{\mu}(\mathfrak{o})G_{\mu,\nu}^{\alt,\lambda}(\mathfrak{o})\\
&=\sum_\mu G_{\mu^{(i_1)},\ldots,\mu^{(i_r)}}^{\mu}(\mathfrak{o})G_{\mu,\nu}^{\alt,\lambda}(\mathfrak{o})\\
&=G_{\mu^{(i_1)},\ldots,\mu^{(i_r)},\nu}^{\alt,\lambda}(\mathfrak{o}).
\end{split}\end{align}
The first and third equality come from \eqref{item: alt and usual polynomial}, and the second comes from \Cref{prop: properties of the Hall polynomial}.

\item Again, we write $G_{\mu^{(1)},\ldots,\mu^{(r)},\nu}^{\alt,\lambda}(\mathfrak{o})=\sum_\mu G_{\mu^{(1)},\ldots,\mu^{(r)}}^{\mu}(\mathfrak{o})G_{\mu,\nu}^{\alt,\lambda}(\mathfrak{o})$. Both $G_{\mu^{(1)},\ldots,\mu^{(r)}}^{\mu}(\mathfrak{o})$ and $G_{\mu,\nu}^{\alt,\lambda}(\mathfrak{o})$ are polynomials in $q$ with integer coefficients, then so is $G_{\mu^{(1)},\ldots,\mu^{(r)},\nu}^{\alt,\lambda}(\mathfrak{o})$.
\end{enumerate}
\end{proof}

\begin{defi}\label{def:alt_hall_module}
Let $H=H(\mathfrak{o})$ be the Hall algebra of $\mathfrak{o}$. Then we define the \emph{alternating Hall module} $H^{\alt}=H^{\alt}(\mathfrak{o})$ to be the $H(\mathfrak{o})$-module with $\Z$-basis $\{u_\lambda^{\alt}\}$, where the indices $\lambda$ are partitions. Define the multiplication of $H(\mathfrak{o})$ over $H^{\alt}(\mathfrak{o})$ by
$$u_\mu u_\nu^{\alt}=\sum_{\lambda\in\Y} G_{\mu,\nu}^{\alt,\lambda}(\mathfrak{o})u_\lambda^{\alt}=\sum_{\lambda\in\Y}  g_{\mu,\nu}^{\alt,\lambda}(q)u_\lambda^{\alt},\quad\forall\mu,\nu\in\Y.$$
\end{defi}

The multiplication action in \Cref{def:alt_hall_module} defines a module structure because the coefficient of $u_{\l}^{\alt}$ in either $(u_{\mu^{(1)}}u_{\mu^{(2)}})u_{\nu}^{\alt}$ or $u_{\mu^{(1)}}(u_{\mu^{(2)}}u_{\nu}^{\alt})$ is $G_{\mu^{(1)},\mu^{(2)},\nu}^{\alt,\l}$. $H^{\alt}(\mathfrak{o})$ has a canonical graded module structure
$$H^{\alt}(\mathfrak{o})=\bigoplus_{k\ge 0}H^{k,\alt}(\mathfrak{o}),$$
where $H^{k,\alt}(\mathfrak{o})$ is the subgroup with free $\Z$-basis $\{u^{\alt}_\lambda\mid\lambda\in\Y, |\lambda|=k\}$.

The alternating Hall module we defined in this section and the Hecke structure studied in \cite[Section 3]{shen2024non} lead to the following theorem, which is the precise version of \eqref{eq:commutative diagrams}.

\begin{theorem}\label{thm: alt_commutative diagram}
The diagram 
\begin{equation}\label{eq:commutative diagram_alternating}\xymatrix@C=0.5em{
H=H(\mathfrak{o})\ar[d]^{(2)} &\overset{(1)}{\scalebox{2}[1]{\(\curvearrowright\)}} & H^{\alt}=H^{\alt}(\mathfrak{o}) \ar[d]^{(3)}\\
H(G^+,K) \ar[d]^{(5)}&\overset{(4)}{\scalebox{2}[1]{\(\curvearrowright\)}} & H(G^{+{\alt}},K)\ar[d]^{(6)}\\
\Lambda_n\otimes_\Z\Q&\overset{(7)}{\scalebox{2}[1]{\(\curvearrowright\)}} &\Lambda_n\otimes_\Z\Q
}\end{equation}
commutes, where the objects and maps are as follows. The three rings in the left column were defined in \Cref{def:hall_alg} (top), \Cref{defi: Hecke algebra} (middle), and \Cref{defi: symmetric functions} (bottom). The curved arrows denote module structures, and the modules in the right column are as in \Cref{def:alt_hall_module}(top), \Cref{defi: alt_integer value Hecke module} (middle), and the trivial module structure of $\Lambda_n \ot_\Z \Q$ over itself (bottom). The vertical maps are defined on generators by

(2) $u_\mu\mapsto c_\mu$,

(3) $u_\nu^{\alt}\mapsto c_\nu^{\alt}$,

(5) $\psi_n:c_\mu\mapsto q^{-n(\mu)+|\mu|}P_\mu(x_1,x_1q^{-1},\ldots,x_n,x_nq^{-1};q^{-1})$,

(6) $\psi_n^{\alt}:c_\nu^{\alt}\mapsto q^{-2n(\nu)}P_\nu(x_1,\ldots,x_n;q^{-2})$,

\noindent and extended by linearity.
\end{theorem}

\begin{proof}
The upper half of the diagram is commutative due to the following two observations:
\begin{enumerate}
\item \Cref{thm: alt_Hall and Hecke}, which proves that the alternating Hall module and the alternating Hecke module share the same structure coefficients $g_{\mu,\nu}^{\alt,\lambda}(q)$, and
\item (2.6) of \cite[Chapter V]{Macdonald}, which proves that the Hall algebra and Hecke ring also have the same structure coefficients.
\end{enumerate}
For the lower half diagram, \cite[Theorem 3.3]{shen2024non} proved that the mappings (here $\langle\cdot,\cdot\rangle$ is the canonical inner product, $\rho_{2n}=\frac{1}{2}(2n-1,2n-3,\ldots,1-2n),\rho_n=\frac{1}{2}(n-1,n-3,\ldots,1-n)$)
\begin{align}\label{eq: alt_Hironaka}
\begin{split}c_\mu&\mapsto q^{\langle \mu,\rho_{2n}\rangle}P_\mu(x_1q^{\frac{1}{2}},x_1q^{-\frac{1}{2}},\ldots,x_nq^{\frac{1}{2}},x_nq^{-\frac{1}{2}};q^{-1})\\
c_\nu^{\alt}&\mapsto q^{2\langle\nu,\rho_n\rangle}P_\nu(x_1,\ldots,x_n;q^{-2})
\end{split}
\end{align}
extended by linearity give a homomorphism from the $H(G^+,K)$-module $H(G^{+{\alt}},K)$ to $\L_n\otimes_\Z\Q$, viewed as a module over itself. Therefore, a change of variables $x_i\mapsto q^{1-n}x_i$ of the homomorphisms in \eqref{eq: alt_Hironaka} provides the mappings (5) and (6) in \eqref{eq:commutative diagram_alternating}.
\end{proof}

\begin{proof}[Proof of \Cref{thm: Hall and Symmetric functions}, alternating case]
For $\phi$ and $\phi^{\alt}$ as defined in \Cref{thm: Hall and Symmetric functions}, we must first show that
\begin{equation}
\phi(u_\mu)\phi^{\alt}(u_\nu^{\alt}) = \phi^{\alt}(u_\mu * u_\nu^{\alt}).    
\end{equation}
Both sides are elements of $\La_\Q$, so it suffices to check that the projections of both sides to $\La_n \ot_\Z \Q$ are equal for every $n$. But since the projections are homomorphisms, this is exactly the statement
\begin{equation}
\phi_n(u_\mu)\phi_n^{\alt}(u_\nu^{\alt}) = \phi_n^{\alt}(u_\mu * u_\nu^{\alt})    
\end{equation}
which we verified in \Cref{thm: alt_commutative diagram}.

% For any $n\ge 1$, the commutative diagram \eqref{eq:commutative diagram_alternating} provides a homomorphism $(\phi_n,\phi_n^{\alt})$ from $H(\mathfrak{o})$-graded module $H^{\alt}(\mathfrak{o})$ to $\L_n\otimes_\Z\Q$. Explicitly, we have
% \begin{align}
% \begin{split}
% \phi_n:H(\mathfrak{o})&\rightarrow \L_n\otimes_\Z\Q\\
% u_\mu&\mapsto q^{-n(\mu)+|\mu|}P_\mu(x_1,x_1q^{-1},\ldots,x_n,x_nq^{-1};q^{-1})
% \end{split}
% \end{align}
% and
% \begin{align}
% \begin{split}
% \phi^{\alt}_n:H^{\alt}(\mathfrak{o})&\rightarrow \L_n\otimes_\Z\Q\\
% u^{\alt}_\nu&\mapsto q^{-2n(\nu)}P_\nu(x_1,\ldots,x_n;q^{-2}).
% \end{split}
% \end{align}
% For any $m\ge n$, let $\rho_{m,n}:\Lambda_m\rightarrow\Lambda_n$ be the same as in \Cref{defi: symmetric functions}. Then it is clear that $\phi_n=\rho_{m,n}\circ\phi_m,\phi^{\alt}_n=\rho_{m,n}\circ\phi^{\alt}_m$. Recall from \Cref{rmk: inverse limit of graded} that $\Lambda_\Q$ is the inverse limit $(\Lambda_n\otimes_\Z\Q)_{n\ge 1}$, viewed as graded modules of itself. Due to the universal mapping property, the maps $\phi_n$ and $\phi_n^{\alt}$ uniquely lift to maps 
% \begin{equation}
%     \phi: H(\mf{o}) \to \Lambda_\Q 
% \end{equation}
% and 
% \begin{equation}
%     \phi^{\alt}: H^{\alt}(\mf{o}) \to \Lambda_\Q,
% \end{equation}
% which are given by \eqref{eq: alt_image of ring} and \eqref{eq: alt_image of module} respectively. 

Now we must prove the maps $(\phi,\phi^{\alt})$ we construct indeed yield an isomorphism. First, we show that $\phi$ provides a ring isomorphism from $H(\mathfrak{o})\otimes_\Z\Q$ to $\L_\Q$. \Cref{cor: Hall and symmetric functions} already gives an isomorphism $\psi$ (taking tensor product with $\Q$) between these two rings, so we only need to show that the $\Q$-linear mappings 
\begin{multline}
\phi\circ\psi^{-1}:q^{-n(\mu)}P_\mu(x_1,x_2,\ldots;q^{-1}) \mapsto \\q^{-n(\mu)+|\mu|}P_\mu(x_1,x_1q^{-1},x_2,x_2q^{-1},\ldots;q^{-1})=q^{-n(\mu)}P_\mu(x_1,x_1q,x_2,x_2q,\ldots;q^{-1})
\end{multline}
actually gives a ring automorphism of $\L_\Q$. In fact, the map could also be expressed in the form $\phi\circ\psi^{-1}(p_n)=(1+q^n)p_n$, where $p_n$ is the $n$th power sum. Since the power sums $p_n$ generate $\L_\Q$, $\phi\circ\psi^{-1}$ must be an automorphism with inverse given by $p_n\mapsto(1+q^n)^{-1}p_n$. 

On the other hand, since the $P_\nu$ form a basis of the vector space $\L_\Q$, the map $\phi^{\alt}$ is both injective and surjective, and hence must be an isomorphism. This completes the proof. 
\end{proof}

\begin{corollary}
$H^{\alt}(\mathfrak{o})\otimes_\Z\Q$ is generated by the set $\{u_{(n)}\mid n\ge 0\}\subset H(\mathfrak{o})$ over $u_0^{\alt}$. 
\end{corollary}
\begin{proof}
    This is clear to see when we apply the isomorphism $(\phi,\phi^{\alt})$, which maps every $u_{(n)}$ to $(1+q^n)P_{(n)}$, $u_0^{\alt}$ to $1$ and view the module structure from the symmetric function point of view.
\end{proof}

\begin{rmk}
It is worth mentioning that $H^{\alt}(\mathfrak{o})$ cannot be generated by $H(\mathfrak{o})$ over $u_0^{\alt}$ if we do not first tensor with $\Q$. In fact, the module generated in this way does not contain the element $u_{(1,1)}^{\alt}$.
\end{rmk}

\section{The Hermitian Hall Module}\label{sec:hermitian}

\subsection{Hermitian module background}

As an analogy of the alternating case, we construct a module $H^{\her}=H^{\her}(\mathfrak{o})$ over the Hall algebra $H(\mathfrak{o})$, which shall be called the \emph{Hermitian Hall module}. 

Let $M$ be a finite $\mathfrak{o}$-module, equipped with a bilinear Hermitian pairing $\langle\cdot,\cdot\rangle:M\times M\rightarrow F/\mathfrak{o}$. Here recall that Hermitian means $\langle x,y \rangle =  \langle y, x \rangle^*$. For any $\mathfrak{o}$-submodule $N$ of $M$, we write $\langle\cdot,\cdot\rangle\mid_N$ to denote the restriction of the Hermitian pairing to $N$. Let $\hat M=\text{Hom}(M,F/\mathfrak{o})$ equipped with the canonical $\mathfrak{o}$-module structure. The following properties are equivalent for the Hermitian pairing:
\begin{enumerate}
\item The map $M\rightarrow\hat M, x\mapsto\langle x,\cdot\rangle$ gives an isomorphism from $M$ to $\hat M$.\label{item: her_iso}

\item If $x\in M$ satisfies $\langle x,y\rangle=0$ for every $y\in M$, then $x=0$. \label{item: her_x=0}
\end{enumerate}

\begin{defi}
If the above properties \eqref{item: her_iso} and \eqref{item: her_x=0} are satisfied, we say the pairing is \emph{regular}, and $M$ has a \emph{Hermitian $\mathfrak{o}$-module} structure. 
\end{defi}

Let $(M,\langle\;,\;\rangle),(M',\langle\;,\;\rangle')$ be two Hermitian $\mathfrak{o}$-modules. We say these two $\mathfrak{o}$-modules are \emph{isomorphic}, if there exists an $\mathfrak{o}$-module isomorphism $\varphi:M\rightarrow M'$ that preserves the pairing, i.e., for all $x,y\in M$,
$$\langle x,y\rangle=\langle\varphi(x),\varphi(y)\rangle'.$$
\begin{theorem}
Every Hermitian $\mathfrak{o}$-module $M$ is isomorphic to one of the form
\begin{equation}
M\cong\mathfrak{o}/\mathfrak{p}^{\lambda_1}\oplus\cdots\oplus\mathfrak{o}/\mathfrak{p}^{\lambda_n}
\end{equation}
for some $\lambda_1\ge\ldots\ge\lambda_n\ge 1$ with the pairing $\langle\cdot,\cdot\rangle$ given by $\langle x,y\rangle=x^*\pi_\lambda^{-1}y$ for any $x=(x_1,\ldots,x_n),y=(y_1,\ldots,y_n)\in M$, where
\begin{equation}
\pi_\lambda=\diag(\pi^{\lambda_1},\ldots,\pi^{\lambda_n}).
\end{equation}
\end{theorem}
\begin{proof}
See \cite[Section 7]{Jacobowitz}.
\end{proof}

In this case, we say $M$ is 
a \emph{Hermitian module} of type $\lambda=(\lambda_1,\ldots,\lambda_n)$. The isomorphism class of Hermitian $\mathfrak{o}$-modules is totally determined by the partition $\lambda$.

For any $\mathfrak{o}$-submodule $N$ of $M$, let
\begin{equation}\label{eqref: her_orthogonal submodule}
N^\perp=\{x\in M\mid\langle x,y\rangle=0,\quad\forall y\in N\}\end{equation}
be the \emph{orthogonal submodule of $N$ in $M$}. In this case, the $\mathfrak{o}$-module homomorphism
$$M\rightarrow\hat N:=\Hom(N,F/\mathfrak{o}),\quad x\mapsto\langle x,\cdot\rangle$$
is surjective with kernel $N^\perp$. Note that the isomorphism $N\cong\hat N$ always holds. This can be seen by observing that the mapping commutes with direct sums, and therefore, it is sufficient to verify the case when $N$ is cyclic, which is straightforward. Therefore, we obtain the following:

\begin{proposition}\label[proposition]{thm:perp_her_props}
If $M$ is a Hermitian $\mathfrak{o}$-module and $N,N'\subset M$, then
\begin{enumerate}
\item $N\cong M/N^\perp,N^\perp\cong M/N$ as $\mathfrak{o}$-modules.

\item $\#N\cdot\#N^\perp=\#M$.

\item $N^{\perp\perp}=N$.

\item If $N$ is regular (i.e., the restriction of the pairing over $N$ is regular), then $N^\perp$ is also regular, and $M\cong N\oplus N^\perp$.

\item $N\subset N'\Leftrightarrow N'^\perp\subset N^\perp$. \label{item: her_include perp}
\end{enumerate}
\end{proposition}

The proof of the above proposition is similar to the alternating case, and we omit it here.

An $\mathfrak{o}$-module $N\subset M$ is \emph{isotropic} if $N\subset N^\perp$, i.e. $(N,N)=0$. In this case, the Hermitian pairing over $M$ induces a reduced Hermitian pairing over $N^\perp/N$, which gives this $\mathfrak{o}$-module a Hermitian structure. Also, by part \eqref{item: her_include perp} of the above theorem, every $\mathfrak{o}$-submodule $P\subset N$ is also isotropic.

\subsection{The Hermitian Hecke module $H(G^{\her},K)$}

The goal of this subsection is to provide background for the Hermitian Hecke module studied in \cite{shen2024non}, which turns out in \Cref{thm: her_Hall and Hecke} to share the same structure coefficient as the Hermitian Hall module $H^{\her}(\mathfrak{o})$.

Let $G=\GL_n(F)$ be the group of all invertible $n\times n$ matrices over $F$. Also, let $G^+=G\cap M_n(\mathfrak{o})$ be the subsemigroup of $G$ consisting of all matrices $x\in G$ with entries $x_{ij}\in\mathfrak{o}$, and let $K=\GL_n(\mathfrak{o})=G^+\cap(G^+)^{-1}$ so that $K$ consists of all $x\in G$ with entries $x_{ij}\in\mathfrak{o}$ and $\det(x)$ a unit in $\mathfrak{o}$. We define $L(G,K)$ (resp. $L(G^+,K)$) and $H(G,K)$ (resp. $H(G^+,K)$) similarly as the alternating case.

\begin{defi}
For all $x\in G$, we write $x^*$ as the conjugate transpose of $x$, and $x^{-*}=(x^{-1})^*$ as the inverse conjugate transpose of $x$. Let $G^{\her}=G\cap\Her_n(F)=\{x\in G\mid x^*=x\}$ denote the set of invertible Hermitian matrices, and $G^{+{\her}}=G^+\cap\Her_n(\mathfrak{o})$ denote the subset of $G$ consisting of matrices with entries in $\mathfrak{o}$. 
\end{defi}

\begin{defi}
Let $L(G^{\her},K)$ (resp. $L(G^{+{\her}},K)$) denote the space of all complex-valued continuous functions of compact support on $G^{\her}$ (resp. $G^{+\her}$) which are invariant with respect to $K$, i.e., such that
$$f(k^*xk)=f(x)$$
for all $x\in G$ (resp. $x\in G^+$) and $k\in K$. We may regard $L(G^{+{\her}},K)$ as the subset of $L(G^{\her},K)$ consisting of functions which are $0$ outside $G^{+\her}$. 

\end{defi}

We define a multiplication of $L(G,K)$ over the module $L(G^{\her},K)$ as follows: for $g\in L(G,K)$, $f\in L(G^{\her},K)$,
$$(g*f)(x)=\int_Gg(y)f(y^{-1}xy^{-*})dy.$$
Since $g$ and $f$ are compactly supported, the integration is over a compact set. This product satisfies $g_1*(g_2*f)=(g_1*g_2)*f$ and $c_0*f=f$, hence defines a $L(G,K)$ (resp. $L(G^+,K)$)-module structure of $L(G^{\her},K)$ (resp. $L(G^{+\her},K)$). (We use the same notation $*$ for the transition involution of a given matrix, the convolution of the Hecke ring itself,  and the convolution of the Hecke ring over the Hermitian Hecke module, but it is easy to distinguish these three based on context). 

Each function $f\in L(G^{\her},K)$ is constant on the orbits $\{k^*xk|k\in K\}$ in $G^{\her}$. These orbits are compact, open, and mutually disjoint. Since $f$ has compact support, it follows that $f$ takes non-zero values on only finitely many orbits $\{k^*xk|k\in K\}$, and hence can be written as a finite linear combination of their characteristic functions. Therefore, the characteristic functions of these orbits in $G^{\her}$ form a $\C$-basis of $L(G^{\her},K)$. 

%If we vary the definition of the module $L(G^{\her},K)$ (resp. $L(G^{+{\her}},K)$) by requiring the functions to take their values in $\Z$ instead of $\C$, the resulting module is the generalization of the Hecke ring over $G^{\her}$ (resp. $G^{+{\her}}$), and we denote it by $H(G^{\her},K)$ (resp. $H(G^{+{\her}},K)$). 
\begin{defi}\label{defi: her_integer value Hecke module}
Let $H(G^{\her},K)$ (resp. $H(G^{+\her},K)$) denote the subspace of $L(G^{\her},K)$ (resp. $L(G^{+{\her}},K)$) consisting of all integer-valued functions, which has a $H(G,K)$ (resp. $H(G^+,K)$)-module structure as above. We may regard $H(G^{+{\her}},K)$ as a subset of $H(G^{\her},K)$.
\end{defi}
Clearly we have
$$L(G^{\her},K)=H(G^{\her},K)\otimes_\Z\C,L(G^{+{\her}},K)=H(G^{+{\her}},K)\otimes_\Z\C.$$
Consider an orbit $\{k^*xk|k\in K\}$, where $x\in G^{\her}$. By multiplying $x$ by a suitable power of $\pi$ (the generator of $\mathfrak{p}$) we can bring $x$ into $G^{+{\her}}$. As stated on page 567 of \cite{Hironaka_Hermitian}, each orbit $\{k^*xk|k\in K\}$ has a unique representative of the form
$$\pi_\lambda=\diag(\pi^{\lambda_1},\ldots,\pi^{\lambda_n}),$$
where $\lambda_1\ge\lambda_2\ge\ldots\ge\lambda_n$. We have $\lambda_n\ge 0$ if and only if $x\in G^{+{\her}}$.

Let $c^{\her}_\lambda$ denote the characteristic function of the orbit $\{k^*\pi_\lambda k|k\in K\}$. Then we have the $c^{\her}_\lambda$ (resp. the $c^{\her}_\lambda$ such that $\lambda_n\ge 0$) form a $\Z$-basis of $H(G^{\her},K)$ (resp. $H(G^{+{\her}},K)$). The characteristic function $c_0$ of $K$ is the identity element of $H(G,K)$ and $H(G^+,K)$, and this also plays the role of the identity element when multiplying with elements in $H(G^{\her},K)$ and $H(G^{+{\her}},K)$. 

Let $\mu,\nu\in\Y_n$ be partitions. Then, the product $c_\mu *c^{\her}_\nu$ has the form
\begin{equation}\label{eq:structure coefficient_her}c_\mu *c^{\her}_\nu=\sum_{\lambda\in\Y_n} g_{\mu,\nu}^{\her,\lambda}(q)c^{\her}_\lambda,
\end{equation}
where \cite[Corollary 4.4]{shen2024non} proves that $g_{\mu,\nu}^{\her,\lambda}(q)\in\Z[q]$ is a polynomial in $q$, independent of $\mathfrak{o}$. Explicitly, we have
$$g_{\mu,\nu}^{\her,\lambda}(q)=(c_\mu*c_\nu^{\her})(\pi_\lambda)=\int_G c_\mu(y)c_\nu^{\her}(y^{-1}\pi_\lambda y^{-*})dy.$$
Since $c_\mu(y)$ vanishes for $y$ outside $K\pi_\mu K$, the integration is over this orbit, which we shall write as a disjoint union of left cosets, say
\begin{equation}\label{eq: her_disjoint left coset}K\pi_\mu K=\bigsqcup_j y_jK\quad (y_j\in K\pi_\mu).\end{equation}
This gives 
\begin{equation}\label{eq: her_sum of left cosets}g_{\mu,\nu}^{\her,\lambda}(q)=\sum_j\int_{y_jK}c^{\her}_\nu(y^{-1}\pi_\lambda y^{-*})dy=\sum_jc^{\her}_\nu(y_j^{-1}\pi_\lambda y_j^{-*}).\end{equation}

\subsection{The Hermitian Hall polynomial}

The following theorem is the critical step that connects the Hermitian Hall module and the Hermitian Hecke module.

\begin{theorem}\label{thm: her_Hall and Hecke}
Let $\lambda,\mu,\nu\in\Y_n$. The polynomial $g_{\mu,\nu}^{\her,\lambda}(q)$ is equal to $G_{\mu,\nu}^{\her,\lambda}(\mathfrak{o})$, the number of submodules $M'\subset M$ such that $M/M'$ is an $\mathfrak{o}$-module of type $\mu$, $M'^\perp$ is isotropic, and $M'/M'^\perp$ is a Hermitian $\mathfrak{o}$-module of type $\nu$.

\end{theorem}

\begin{proof} 
Consider the lattice $L=\bigoplus_{i=1}^n\mathfrak{o}$ equipped with a Hermitian pairing $L\times L\rightarrow F/\mathfrak{o}: \langle x,y\rangle= x^*\pi_\lambda^{-1}y$, where $x,y\in L$. We view $M$ as the quotient $M=L/\pi_\lambda L$ with bilinear pairing inherited from $L$. It is clear that the left cosets $y_jK$ in \eqref{eq: her_disjoint left coset} are in one-to-one correspondence with the submodules $M'\subset M$ such that $M/M'$ is a finite $\mathfrak{o}$-module with type $\mu$, under the map $y_j K \mapsto M'=y_jL/\pi_\lambda L$. Then we claim the following are equivalent:
\begin{enumerate}[left=0pt]
\item $c_\nu(y_j^{-1}\pi_\lambda y_j^{-*})=1$.\label{item: her_c=1}

\item There exists $k_j\in K$ such that $y_j^{-1}\pi_\lambda y_j^{-*}=k_j\pi_\nu k_j^*$.\label{item: her_existence}

\item $M'\supset M'^\perp$, and $M'/M'^\perp$ is a Hermitian module of type $\nu$. Here $M'^\perp$ is the orthogonal submodule of $M'$ defined in \eqref{eqref: her_orthogonal submodule}. Note that it follows from the previous sentence that $M'^\perp$ is isotropic.
\label{item: her_type}
\end{enumerate}

Now let us prove the above equivalence:

\eqref{item: her_c=1}$\Rightarrow$\eqref{item: her_existence}: Trivial.

\eqref{item: her_existence}$\Rightarrow$\eqref{item: her_type}: In this case, we have $M'=y_jL/\pi_\lambda L$. We claim that $M'^\perp=y_jk_j\pi_\nu k_j^*L/\pi_\lambda L$. An element $y_j a \in M'$ will have
\begin{equation*}
    \langle \pi_\la y_j^{-*} b, y_j a \rangle = b^* y_j^{-1} \pi_\la \pi_\la^{-1} y_j a = b^* a,
\end{equation*}
which is $0$ in $F/\mf{o}$ for all $a,b\in\bigoplus_{i=1}^n\mathfrak{o}$, hence $\pi_\lambda y_j^{-*}L/\pi_\lambda L \subset M'^\perp$.

Conversely, if $b\in\bigoplus_{i=1}^n F$ satisfies $b^*a\in\mathfrak{o}$ for all $a\in \bigoplus_{i=1}^n\mathfrak{o}$, then $b\in\bigoplus_{i=1}^n\mathfrak{o}$, which shows the opposite inclusion. Thus
\begin{equation}
    \label{eq:mperp_herm}
    M'^\perp=\pi_\lambda y_j^{-*}L/\pi_\lambda L=y_jk_j\pi_\nu k_j^*L/\pi_\lambda L,
\end{equation}
where the second equality is by the hypothesis \eqref{item: her_existence}. The right hand side of \eqref{eq:mperp_herm} is contained in $M'=y_jL/\pi_\lambda L$, so $M'\supset M'^\perp$. Finally, the Hermitian pairing restricted to $M'$ is given by 
\begin{equation}
    \langle y_jk_jx,y_jk_jy\rangle= x^* k_j^* y_j^* \pi_\la^{-1} y_j k_j y = x^*\pi_\nu ^{-1}y
\end{equation}
by \eqref{item: her_existence}, which shows that $M'/M'^\perp$ is a Hermitian module of type $\nu$.

\eqref{item: her_type}$\Rightarrow$\eqref{item: her_c=1}: In this case, the type of the Hermitian module $M'/M'^\perp$ implies the Hermitian pairing restricted over $M'$ has type $\nu$, i.e., $c_\nu(y_j^{-1}\pi_\lambda y_j^{-*})=1$. 

Therefore, following the expression in \eqref{eq: her_sum of left cosets}, the number of $y_j$ that satisfies the above equivalent properties \eqref{item: her_c=1}, \eqref{item: her_existence}, and \eqref{item: her_type} is $G_{\mu,\nu}^{\her,\lambda}(\mathfrak{o})$. This ends the proof.
\end{proof}

\Cref{thm: her_Hall and Hecke} motivates us to generalize the coefficient $G_{\mu,\nu}^{\her,\lambda}(\mathfrak{o})=g_{\mu,\nu}^{\her,\lambda}(q)$ which appeared in \eqref{eq:structure coefficient_her}, based on the Hermitian Hall module we defined.

\begin{defi}

Let $\mu^{(1)},\ldots,\mu^{(r)},\lambda,\nu\in\Y_n$. We define $G_{\mu^{(1)},\ldots,\mu^{(r)},\nu}^{\her,\lambda}(\mathfrak{o})$ to be the number of chains of submodules of $M$:
$$M=M_0\supset M_1\supset\ldots\supset M_r$$
such that $M_{i-1}/M_i$ has type $\mu^{(i)}$ as an $\mathfrak{o}$-module, $M_i^\perp$(viewed as the orthogonal submodule of $M_i$ in $M$) is isotropic for $1\le i\le r$, and $M_r/M_r^\perp$ is of type $\nu$.   

\end{defi}

\begin{proposition}
$G_{\mu^{(1)},\ldots,\mu^{(r)},\nu}^{\her,\lambda}(\mathfrak{o})$ has the following properties:
\begin{enumerate}
\item $G_{\mu^{(1)},\ldots,\mu^{(r)},\nu}^{\her,\lambda}(\mathfrak{o})=0$ unless $|\lambda|=2|\mu^{(1)}|+\cdots+2|\mu^{(r)}|+|\nu|$.

\item $G_{\mu^{(1)},\ldots,\mu^{(r)},\nu}^{\her,\lambda}(\mathfrak{o})=\sum_\mu G_{\mu^{(1)},\ldots,\mu^{(r)}}^{\mu}(\mathfrak{o})G_{\mu,\nu}^{\her,\lambda}(\mathfrak{o})$, where the sum covers $\mu\in\Y$ with $l(\mu)\le n$.\label{item: her and usual polynomial}

\item $G_{\mu^{(1)},\ldots,\mu^{(r)},\nu}^{\her,\lambda}(\mathfrak{o})$ does not depend on the order of $\mu^{(1)},\ldots,\mu^{(r)}$, i.e., for any permutation $\{i_1,\ldots,i_r\}=\{1,2,\ldots,r\}$, we have $G_{\mu^{(1)},\ldots,\mu^{(r)},\nu}^{\her,\lambda}(\mathfrak{o})=G_{\mu^{(i_1)},\ldots,\mu^{(i_r)},\nu}^{\her,\lambda}(\mathfrak{o})$.

\item $G_{\mu^{(1)},\ldots,\mu^{(r)},\nu}^{\her,\lambda}(\mathfrak{o})$ is a ``polynomial in $q$", i.e., there exists a polynomial $g_{\mu^{(1)},\ldots,\mu^{(r)},\nu}^{\her,\lambda}(t)\in\Z[t]$, independent of $\mathfrak{o}$, such that $G_{\mu^{(1)},\ldots,\mu^{(r)},\nu}^{\her,\lambda}(\mathfrak{o})=g_{\mu^{(1)},\ldots,\mu^{(r)},\nu}^{\her,\lambda}(q)$. 
\end{enumerate}
\end{proposition}

\begin{proof}
\begin{enumerate}
\item This is because for all $M=M_0\supset M_1\supset\ldots\supset M_r$ such that $M_{i-1}/M_i$ has type $\mu^{(i)}$ as an $\mathfrak{o}$-module, $M_i^\perp$ is isotropic for $1\le i\le r$, and $M_r/M_r^\perp$ is of type $\nu$ as a Hermitian $\mathfrak{o}$-module, we have 
\begin{align*}
|\lambda|=l(M_0)=\sum_{i=1}^{r}(l(M_{i-1}/M_i)+l(M_i^\perp/M_{i-1}^\perp))+l(M_r/M_r^\perp)=2|\mu^{(1)}|+\cdots+2|\mu^{(r)}|+|\nu|,
\end{align*}
where $l(\cdot)$ is the length defined in \eqref{eq:def_length}.

\item To count the number of chains of modules $M=M_0\supset M_1\supset\ldots\supset M_r$ that satisfy our requirement, we may start from finding a submodule $M_r\subset M$ such that $M_r^\perp$ is isotropic, and $M_r/M_r^\perp$ is a Hermitian $\mathfrak{o}$-module of type $\nu$. Then we choose the remaining elements of the chain $M=M_0\supset M_1\supset\ldots\supset M_r$ such that $M_{i-1}/M_i$ has type $\mu^{(i)}$ as an $\mathfrak{o}$-module. Since $M_r^\perp$ is isotropic and $M_i^\perp \subset M_r^\perp$, $M_i^\perp$ is automatically isotropic for all $1\le i\le r-1$ by \Cref{thm:perp_her_props}. Summing over all the types of $M/M_r$ as an $\mathfrak{o}$-module brings the result.
\item We have
\begin{align}\begin{split}
G_{\mu^{(1)},\ldots,\mu^{(r)},\nu}^{\her,\lambda}(\mathfrak{o})&=\sum_\mu G_{\mu^{(1)},\ldots,\mu^{(r)}}^{\mu}(\mathfrak{o})G_{\mu,\nu}^{\her,\lambda}(\mathfrak{o})\\
&=\sum_\mu G_{\mu^{(i_1)},\ldots,\mu^{(i_r)}}^{\mu}(\mathfrak{o})G_{\mu,\nu}^{\her,\lambda}(\mathfrak{o})\\
&=G_{\mu^{(i_1)},\ldots,\mu^{(i_r)},\nu}^{\her,\lambda}(\mathfrak{o}).
\end{split}\end{align}
The first and third equality come from \eqref{item: her and usual polynomial}, and the second comes from \Cref{prop: properties of the Hall polynomial}.

\item Again, we write $G_{\mu^{(1)},\ldots,\mu^{(r)},\nu}^{\her,\lambda}(\mathfrak{o})=\sum_\mu G_{\mu^{(1)},\ldots,\mu^{(r)}}^{\mu}(\mathfrak{o})G_{\mu,\nu}^{\her,\lambda}(\mathfrak{o})$. Both $G_{\mu^{(1)},\ldots,\mu^{(r)}}^{\mu}(\mathfrak{o})$ and $G_{\mu,\nu}^{\her,\lambda}(\mathfrak{o})$ are polynomial in $q$ with integer coefficients, then so is $G_{\mu^{(1)},\ldots,\mu^{(r)},\nu}^{\her,\lambda}(\mathfrak{o})$.
\end{enumerate}
\end{proof}

\begin{defi}\label{def:her_hall_module}
Let $H=H(\mathfrak{o})$ be the Hall algebra of $\mathfrak{o}$. Then we define the \emph{Hermitian Hall module} $H^{\her}=H^{\her}(\mathfrak{o})$ be the $H(\mathfrak{o})$-module with $\Z$-basis $\{u_\lambda^{\her}\}$, where the index $\lambda$ are partitions. Define the action by multiplication of $H(\mathfrak{o})$ on $H^{\her}(\mathfrak{o})$ by
$$u_\mu u_\nu^{\her}=\sum_{\lambda\in\Y} G_{\mu,\nu}^{\her,\lambda}(\mathfrak{o})u_\lambda^{\her}=\sum_{\lambda\in\Y} g_{\mu,\nu}^{\her,\lambda}(q)u_\lambda^{\her},\quad\forall\mu,\nu\in\Y.$$
\end{defi}

The multiplication action in \Cref{def:her_hall_module} defines a module structure because the coefficient of $u_{\l}^{\her}$ in either $(u_{\mu^{(1)}}u_{\mu^{(2)}})u_{\nu}^{\her}$ or $u_{\mu^{(1)}}(u_{\mu^{(2)}}u_{\nu}^{\her})$ is $G_{\mu^{(1)},\mu^{(2)},\nu}^{\her,\l}$. The module $H^{\her}(\mathfrak{o})$ has a natural decomposition
$$H^{\her}(\mathfrak{o})=H_\even^{\her}(\mathfrak{o})\oplus H_\odd^{\her}(\mathfrak{o}),$$
where both $H_\even^{\her}(\mathfrak{o})$ and $H_\odd^{\her}(\mathfrak{o})$ have graded module structures
$$H_\even^{\her}(\mathfrak{o})=\bigoplus_{k\ge 0}H^{2k,\her}(\mathfrak{o}),\quad H_\odd^{\her}(\mathfrak{o})=\bigoplus_{k\ge 0}H^{2k+1,\her}(\mathfrak{o}),$$
and $H^{k,\her}(\mathfrak{o})$ is the subgroup with free $\Z$-basis $\{u^{\her}_\lambda\mid\lambda\in\Y, |\lambda|=k\}$.

The Hermitian Hall module we defined in this section and the Hecke structure studied in \cite[Section 4]{shen2024non} lead to the following theorem, which is the Hermitian version of \eqref{eq:commutative diagrams}.

\begin{theorem}\label{thm: her_commutative diagram}
The diagram 
\begin{equation}\label{eq:commutative diagram_Hermitian}\xymatrix@C=0.5em{
H=H(\mathfrak{o})\ar[d]^{(2)} &\overset{(1)}{\scalebox{2}[1]{\(\curvearrowright\)}} & H^{\her}=H^{\her}(\mathfrak{o}) \ar[d]^{(3)}\\
H(G^+,K) \ar[d]^{(5)}&\overset{(4)}{\scalebox{2}[1]{\(\curvearrowright\)}} & H(G^{+{\her}},K)\ar[d]^{(6)}\\
\Lambda_n^\sq\otimes_\Z\Q&\overset{(7)}{\scalebox{2}[1]{\(\curvearrowright\)}} &\Lambda_n\otimes_\Z\Q
}\end{equation}
commutes, where the objects and maps are as follows. The three rings in the left column were defined in \Cref{def:hall_alg} (top), \Cref{defi: Hecke algebra} (middle), and \Cref{def:Lambda_sq} (bottom). The curved arrows denote module structures, and the modules in the right column are as in \Cref{def:her_hall_module} (top), \Cref{defi: her_integer value Hecke module} (middle), and the trivial module structure of $\Lambda_n^{\sq} \ot_\Z \Q$ over $\Lambda_n \ot_\Z \Q$ (bottom). The vertical maps are defined on generators by

(2) $u_\mu\mapsto c_\mu$,

(3) $u_\nu^{\her}\mapsto c_\nu^{\her}$,

(5) $\psi_n:c_\mu\mapsto q^{-2n(\mu)}P_\mu(x_1^2,\ldots,x_n^2;q^{-2})$,

(6) $\psi^{\her}_n:c_\nu^{\her}\mapsto (-1)^{n(\nu)}q^{-n(\nu)}P_\nu(x_1,\ldots,x_n;-q^{-1})$.

\noindent and extended by linearity.

\end{theorem}

\begin{proof}
The upper half of the diagram is commutative due to the following two observations:
\begin{enumerate}
\item \Cref{thm: her_Hall and Hecke}, which proves that the Hermitian Hall module and the Hermitian Hecke module share the same structure coefficients $g_{\mu,\nu}^{\her,\lambda}(q)$; and
\item (2.6) of \cite[Chapter V]{Macdonald}, which proves that the Hall algebra and Hecke ring also have the same structure coefficients.
\end{enumerate}
For the lower half diagram, \cite[Theorem 4.3]{shen2024non} proved that the mappings (here $\langle\cdot,\cdot\rangle$ is the canonical inner product, $\rho_n=\frac{1}{2}(n-1,n-3,\ldots,1-n)$)
\begin{align}\label{eq: her_Hironaka}
\begin{split}c_\mu&\mapsto q^{2\langle \mu,\rho_{n}\rangle}P_\mu(x_1^2,\ldots,x_n^2;q^{-2})\\
c_\nu^{\her}&\mapsto (-1)^{n(\nu)+|\nu|}q^{\langle\nu,\rho_n\rangle}P_\nu(x_1,\ldots,x_n;-q^{-1})
\end{split}
\end{align}
extended by linearity give a homomorphism from the $H(G^+,K)$-module $H(G^{+{\her}},K)$ to $\L_n\otimes_\Z\Q$, viewed as a module over itself. Therefore, a change of variables $x_i\mapsto -q^{\frac{1-n}{2}}x_i$ of the homomorphisms in \eqref{eq: her_Hironaka} provides the mappings (5) and (6) in \eqref{eq:commutative diagram_Hermitian}.
\end{proof}

\begin{proof}[Proof of \Cref{thm: Hall and Symmetric functions}, Hermitian case]
For $\phi$ and $\phi^{\her}$ as defined in \Cref{thm: Hall and Symmetric functions}, we must first show that
\begin{equation}
\phi(u_\mu)\phi^{\her}(u_\nu^{\her}) = \phi^{\her}(u_\mu * u_\nu^{\her}).    
\end{equation}
Both sides are elements of $\La_\Q$, so it suffices to check that the projections of both sides to $\La_n \ot_\Z \Q$ are equal for every $n$. But since the projections are homomorphisms, this is exactly the statement
\begin{equation}
\phi_n(u_\mu)\phi_n^{\her}(u_\nu^{\her}) = \phi_n^{\her}(u_\mu * u_\nu^{\her})    
\end{equation}
which we verified in \Cref{thm: her_commutative diagram}.

Now we must prove the map $(\phi,\phi^{\her})$ we construct is indeed an isomorphism. First, we show that $\phi$ provides a ring isomorphism from $H(\mathfrak{o})\otimes_\Z\Q$ to $\L_\Q^{\sq}$. Since the $P_\mu(x_1,x_2,\ldots;q^{-2})$ form a basis of the vector space $\L_\Q$, the $P_\mu(x_1^2,x_2^2,\ldots;q^{-2})$ must form a basis of the vector space $\L_\Q^{\sq}$. Therefore, the map $\phi$ is an isomorphism because it is both injective and surjective. On the other hand, since the $P_\nu(x_1,x_2,\ldots;-q^{-1})$ form a basis of the vector space $\L_\Q$, the map $\phi^{\her}$ is both injective and surjective, and hence must be an isomorphism. This completes the proof.
\end{proof}

\begin{rmk}
$H^{\her}(\mathfrak{o})\otimes_\Z\Q$ cannot be generated by the scalar multiplication of $H(\mathfrak{o})$ over $u_0^{\her},u_{(1)}^{\her}$. When we apply the isomorphism $(\phi,\phi^{\her})$, it is clear that $P_{(1,1)}(x_1,x_2,\ldots;-q^{-1})=\sum_{i<j}x_ix_j$ is not included in $\L_\Q^{\sq}$, hence we cannot get $u_{(1,1)}^{\her}$ by the scalar multiplication of $H(\mathfrak{o})$ over $u_0^{\her}$. Likewise, we cannot get $u_{(1,1,1)}^{\her}$ by the scalar multiplication of $H(\mathfrak{o})$ over $u_{(1)}^{\her}$. 
\end{rmk}

\section{The $u$-probability}\label{sec:u}

As we have seen in \Cref{subsec: Random modules with pairings}, the formula for $u$-probability involves the number of automorphisms, while for the alternating and Hermitian case, the automorphisms are required to preserve the pairing. The following theorem shows that the $u$-measures for the usual, alternating, and Hermitian cases are all Hall-Littlewood measures. It also explains the motivation for us to study these $u$-probabilities: they come from the distribution of cokernels of the usual, alternating, and Hermitian random matrices. 

\begin{theorem}\label{thm:u-measure_formulas}
Let $t=q^{-1}$ and $\mathbf{P}, \mathbf{P}^{\alt},\mathbf{P}^{\her}$ be the probability measures on partitions defined in \eqref{eq:P}, \eqref{eq:P_alt}, and \eqref{eq:P_her} respectively. Then we have explicit formulas
\begin{enumerate}
\item (No-pairing case) $$\mathbf{P}(u;\lambda)=\frac{q^{-u|\l|}\prod_{j\ge 1}(1-q^{-u-j})}{q^{|\lambda|+2n(\lambda)}\prod_{i\ge 1}(q^{-1};q^{-1})_{m_i(\lambda)}}=\frac{P_\lambda(t^{1+u},t^{2+u},\ldots;t)Q_\lambda(1,t,\ldots;t)}{\Pi_t(t^{1+u},t^{2+u},\ldots;1,t,\ldots)}.$$ 
%In particular, as proved in \cite[Corollary 1.4]{van2020limits}, if $u=k\ge 0$ is an integer, the above becomes the limit of the distribution of $\SN(A_{n\times(n+k)})$ when $n$ goes to infinity, where $A_{n\times(n+k)}$ is a random $n\times (n+k)$ matrix which entries are i.i.d. Haar distributed in $\mathfrak{o}$;
\item (Alternating case)$$\mathbf{P}^{\alt}(u;\lambda)=\frac{q^{(2-2u)|\lambda|}\prod_{j\ge 1}(1-q^{-2u-2j+1})}{q^{4n(\lambda)+3|\lambda|}\prod_{i\ge 1}(q^{-2};q^{-2})_{m_i(\lambda)}}=\frac{P_\lambda(t^{1+2u},t^{3+2u},\ldots;t^2)Q_\lambda(1,t^2,\ldots;t^2)}{\prod_{t^2}(t^{1+2u},t^{3+2u},\ldots;1,t^2,\ldots)}.$$
%In particular, as proved in \cite[Theorem 3.9]{bhargava2013modeling}, if $u=0$, the above becomes the limit $\mathbf{P}_{2\infty}^{\alt}(\lambda)=\lim_{n\rightarrow\infty}\mathbf{P}_{2n}^{\alt}(\lambda)$; Furthermore, when $u=1$, \cite{shen2024non} shows that the above becomes the limit $\mathbf{P}_{2\infty+1}^{\alt}(\lambda)=\lim_{n\rightarrow\infty}\mathbf{P}_{2n+1}^{\alt}(\lambda)$;
\item (Hermitian case)
$$\mathbf{P}^{\her}(u;\lambda)=\frac{q^{-u|\lambda|}\prod_{j\ge 1}(1-(-1)^{j-1}q^{-u-j})}{q^{|\lambda|+2n(\lambda)}\prod_{i\ge 1}(-q^{-1};-q^{-1})_{m_i(\lambda)}}=\frac{P_\lambda(t^{1+u},-t^{2+u},\ldots;-t)Q_\lambda(1,-t,\ldots;-t)}{\prod_{-t}(t^{1+u},-t^{2+u},\ldots;1,-t,\ldots)}.$$
%In particular, as proved in \cite[(1.8)]{shen2024non}, if $u=0$, the above becomes the limit $\mathbf{P}_\infty^{\her}(\lambda)=\lim_{n\rightarrow\infty}\mathbf{P}_n^{\her}(\lambda)$.
\end{enumerate}
\end{theorem}

\begin{proof}
For the usual case, (1.6) of \cite[Chapter II]{Macdonald} gives the number of automorphisms $\#\Aut(M_\lambda)=q^{2n(\lambda)+|\lambda|}\prod_{i\ge 1}(q^{-1};q^{-1})_{m_i(\lambda)}$. For the alternating case, \cite[Theorem 3]{delaunay} gives the number of pairing-preserving automorphisms $\#\Aut^A(M_\lambda)=p^{4n(\lambda)+3|\lambda|}\prod_{i\ge 1}(p^{-2};p^{-2})_{m_i(\lambda)}$ for abelian $p$-groups, but it is clear one can generalize the proof to all alternating $\mathfrak{o}$-modules by replacing $p$ by the $q$, i.e., order of the residue field. For the Hermitian case, we will later prove the explicit form of the number of pairing-preserving automorphisms $\#\Aut^H(M_\lambda)=q^{|\lambda|+2n(\lambda)}\prod_{i\ge 1}(-q^{-1};-q^{-1})_{m_i(\lambda)}$ in \Cref{thm: number of automorphism_Hermitian}. Applying the principal specialization formula in \eqref{eq: principal specialization formula for P} and \eqref{eq: principal specialization formula for Q} gives the proof. 
\end{proof}

%We note that in the alternating case, the coincidence between the explicit formulas of \cite{bhargava2013modeling} and the Hall-Littlewood measure was also noticed by Fulman-Kaplan \cite{fulman2018random}. In the no-pairing case this was noticed earlier by Fulman \cite{fulman-CL} following an intermediate observation of Lengler \cite{lengler_thesis}.

\begin{rmk}\label{rmk:u_meas_and_cokernels}
The $u$-probabilities, in fact, gives the limiting distribution of cokernels of random matrices. For any matrix $A \in \Mat_{n \times (n+k)}(\mf{o})$, its cokernel is the $\mf{o}$-module
\begin{equation}
    \cok(A) := \mf{o}^n / A \mf{o}^{n+k}.
\end{equation}
If $A$ is alternating (resp. Hermitian), $\cok(A)$ inherits the structure of an alternating (resp. Hermitian) pairing, given by $\langle x,y\rangle = x^TA_{2n}^{-1}y$ (resp. $\langle x,y\rangle = x^*A_n^{-1}y$).

Explicitly,
\begin{enumerate}
\item (No-pairing case) For all integer $k\ge 0$, $\mathbf{P}(k;\lambda)$ gives the limit of the distribution of type of $\cok(A_{n \times (n+k)})$ when $n$ goes to infinity, where $A_{n\times(n+k)}$ is a random $n\times (n+k)$ matrix which entries are i.i.d. Haar distributed in $\mathfrak{o}$ with $q=\#(\mathfrak{o}/\mathfrak{p})$ as the order of the residue field. %As proved in \cite[Corollary 1.4]{van2020limits}, f
\item (Alternating case) $\mathbf{P}^{\alt}(0;\lambda)$ gives the limit of the distribution of type of alternating $\mathfrak{o}$-module $\cok(A_{2n})$ when $n$ goes to infinity, where $A_{2n}$ is a random $2n\times 2n$ alternating matrix which above-diagonal entries are i.i.d. Haar distributed in $\mathfrak{o}$. %As proved in \cite[(1.7)]{shen2024non}, 
\item (Hermitian case) $\mathbf{P}^{\her}(0;\lambda)$ gives the limit of the distribution of Hermitian $\mathfrak{o}$-module $\cok(A_n)$ when $n$ goes to infinity, where $A_n$ is a random $n\times n$ Hermitian matrix which above-diagonal entries are i.i.d. Haar distributed in $\mathfrak{o}$, and the diagonal entries are i.i.d. Haar distributed in $\mathfrak{o}\cap F_0$.
\end{enumerate}
The no-pairing and alternating cases are noted in \cite{fulman2018random}, while the Hermitian case is shown in \cite[(1.8)]{shen2024non}. Similar Hall-Littlewood formulas exist at finite $n$ before taking the limit, see \cite[Corollary 1.4]{van2020limits}, \cite{fulman2018random} and \cite{shen2024non} for the no-pairing, alternating, and Hermitian cases respectively.
\end{rmk}

\begin{rmk}
    \label{rmk:u_probs_in_finite_field}
    The $u$-probabilities also govern limits of conjugacy classes of random matrices over finite fields. Up to cosmetic changes of the $u$ variable, the right hand side of each formula in \Cref{thm:u-measure_formulas} is the same with different values of the Hall-Littlewood parameter; this family of measures was studied in many works of Fulman such as \cite{fulman1999probabilistic}. The versions in parts (1) and (3) of \Cref{thm:u-measure_formulas} describe the $N \to \infty$ distribution of the sizes of unipotent Jordan blocks of large random elements of $\GL_N(\F_q)$ and $U_N(\F_q)$ respectively, see \cite{fulman_main} and \cite[Section 4.4]{fulman_thesis}.
\end{rmk}

To compute the number of automorphisms for the Hermitian $\mathfrak{o}$-module $M_\l$, which is needed for \Cref{thm:u-measure_formulas}, we begin with the following technical lemma.

\begin{lemma}\label[lemma]{lem: ways for the first pick}
Let $M_\lambda$ be a Hermitian $\mathfrak{o}$-module of type $\lambda=(\lambda_1,\ldots,\lambda_n)$, where $\lambda_1\ge\ldots\ge\lambda_n\ge 1$. Denote $m=m_{\l_1}(\l)=\#\{i\mid\l_i=\l_1\}$. Then we have
$$\#\{x=(a_1,\ldots,a_n)\in M_\l\mid\langle x,x\rangle=\pi^{-\l_1}\}=q^{2|\lambda|-\lambda_1}(1-(-q)^{-m}).$$
\end{lemma}

\begin{proof}
We first solve the related problem of counting the number of sequences $(a_1,\ldots,a_n)$ such that $a_i\in\mathfrak{o}/\mathfrak{p}^{\lambda_i}$ for all $1\le i\le n$, and
$$\sum_{i=1}^n\pi^{-\lambda_i}\Nm(a_i)=\pi^{-\lambda_1},$$
where $\Nm(a) := a a^*$ is the norm map. While direct computation is difficult, we instead start from the number of sequences $(a_1,\ldots,a_n)$ such that 
\begin{equation}\label{eq:val_lambda1}
    v\left(\sum_{i=1}^n\pi^{-\lambda_i}\Nm(a_i)\right)=-\lambda_1.
\end{equation}
This is equivalent to the condition
\begin{equation}\label{eq:val_0}
    v\left(\sum_{i=1}^m \Nm(a_i)\right)=0,
\end{equation}
which is necessary for \eqref{eq:val_lambda1} to hold since $v(\pi^{-\lambda_i}\Nm(a_i)) \geq -\lambda_i$, and conversely is sufficient since if \eqref{eq:val_0} holds then the other summands in \eqref{eq:val_lambda1} with $\lambda_i < \lambda_1$ do not affect the valuation of the sum.

Therefore, it is natural for us to split the parts in the partition into two subsets $\{\lambda_j\mid\lambda_j=\lambda_1\},\{\lambda_j\mid\lambda_j<\lambda_1\}$ and study them separately. The parts $\lambda_j<\lambda_1$ have no effect on the valuation $v(\cdot)$ of the sum as explained above, and there are $q^{2\sum_{\lambda_j<\lambda_1}\lambda_j}=q^{2|\lambda|-2m\lambda_1}$ ways for us to choose those elements. For the parts $\lambda_j=\lambda_1$, we count the number of sequences $(a_1,\ldots,a_m)$ such that $a_i\in\mathfrak{o}/\mathfrak{p}^{\lambda_1}$ for all $1\le i\le m$, and $v(\sum_{i=1}^m \Nm(a_i))=0$. 

We claim that there are $$q^{2m\lambda_1}(1-(-q)^{-m})(1-q^{-1})$$ 
such sequences $(a_1,\ldots,a_m)$ in total. This claim is trivial when $m=0$. Suppose we have already verified this claim for some $m-1$, then we consider the $m$ case. There are two possible ways that a tuple $(a_1,\ldots,a_m)$ can satisfy the conditions:

\begin{enumerate}
\item $v(\sum_{i=1}^{m-1} \Nm(a_i))>0$, and $v(\Nm(a_m))=0$. By the induction hypothesis, there are
\begin{equation}
    q^{2(m-1)\lambda_1}(1-(1-(-q)^{1-m})(1-q^{-1}))
\end{equation}
choices for $a_1,\ldots,a_{m-1}$ in this case. Now suppose that we have already picked $a_1,\ldots,a_{m-1}$, and need to find $a_m$ such that $v(\Nm(a_m))=0$. Because this is equivalent to $v(a_m) = 0$, there are
\begin{equation}
    q^{2\l_1}(1-q^{-2}) = \#\mf{o}^\times/\mf{p}^{\la_1}
\end{equation}
choices. Hence, there are 
\begin{equation}
    q^{2(m-1)\lambda_1}(1-(1-(-q)^{1-m})(1-q^{-1}))\cdot q^{2\l_1}(1-q^{-2})
\end{equation}
sequences in this case in total. 

\item $v(\sum_{i=1}^{m-1} \Nm(a_i))=0$. Suppose that we already picked $a_1,\ldots,a_{m-1}$. By the induction hypothesis, there are $q^{2(m-1)\lambda_1}(1-(-q)^{1-m})(1-q^{-1})$ choices. Notice that the norm map $\Nm$ sends $\mathfrak{o}^\times/\mathfrak{p}$ to $(\mathfrak{o}\cap F_0)^\times/(\mathfrak{p}\cap F_0)$, where we abusively use quotient notation to denote the image of subsets of $\mf{o}$ under the quotient map, and similarly for $\mf{o} \cap F_0$. The norm map on $\mathfrak{o}^\times/\mathfrak{p}$ is surjective, and it is $(q+1)$ to one since 
\begin{equation}
    \frac{\#\mathfrak{o}^\times/\mathfrak{p}}{\#(\mathfrak{o}\cap F_0)^\times/(\mathfrak{p}\cap F_0)} = \frac{q^2-1}{q-1} = q+1.
\end{equation}
Therefore, the norm also gives a surjective map from $\mathfrak{o}^\times/\mathfrak{p}^{\l_1}$ to $(\mathfrak{o}\cap F_0)^\times/(\mathfrak{p}\cap F_0)^{\l_1}$, which $q^{\l_1}(1+q^{-1})$ to one.

Now, let us consider the choices of $a_m$ such that $v(\sum_{i=1}^{m} \Nm(a_i))>0$. This is equivalent to
\begin{equation}\label{eq:norm_inside}
    \Nm(a_m)\in \left(\mathfrak{p}\cap F_0 -\sum_{i=1}^{m-1} \Nm(a_i)\right)/(\mathfrak{p}\cap F_0)^{\l_1} ,
\end{equation}
where the right hand side is the subset of $(\mathfrak{o}\cap F_0)^\times/(\mathfrak{p}\cap F_0)^{\l_1}$ given by the image of the coset $\mathfrak{p}\cap F_0 -\sum_{i=1}^{m-1} \Nm(a_i)$. This set on the right hand side of \eqref{eq:norm_inside} has size $q^{\l_1-1}$. Hence, there are
\begin{equation}
    q^{\l_1}(1+q^{-1})\cdot q^{\l_1-1}=q^{2\l_1}(q^{-1}+q^{-2})
\end{equation}
choices of $a_m$ such that $v(\sum_{i=1}^{m} \Nm(a_i))>0$. Conversely, the number of choices such that $v(\sum_{i=1}^{m} \Nm(a_i))=0$ is 
\begin{equation}
    q^{2\l_1}-q^{2\l_1}(q^{-1}+q^{-2})=q^{2\l_1}(1-q^{-1}-q^{-2}).
\end{equation}
To sum up, there are
\begin{equation}
    q^{2(m-1)\lambda_1}(1-(-q)^{1-m})(1-q^{-1})\cdot q^{2\l_1}(1-q^{-1}-q^{-2})
\end{equation}
sequences $(a_1,\ldots,a_m)$ in this case in total.
\end{enumerate}
The sum of the above two cases equals \begin{multline*}
q^{2m\lambda_1}((1-(1-(-q)^{1-m})(1-q^{-1}))\cdot (1-q^{-2})+(1-(-q)^{1-m})(1-q^{-1})\cdot(1-q^{-1}-q^{-2}))\\
=q^{2m\lambda_1}(1-(-q)^{-m})(1-q^{-1}).
\end{multline*}
Therefore, by induction, we verify the claim. All in all, we know that there are
$$q^{2m\lambda_1}(1-(-q)^{-m})(1-q^{-1})q^{2|\lambda|-2m\lambda_1}=q^{2|\lambda|}(1-(-q)^{-m})(1-q^{-1})$$
ways to choose the sequence $(a_1,\ldots,a_n)$ such that $v(\sum_{i=1}^n\pi^{-\lambda_i}\Nm(a_i))=-\lambda_1$. 

The sum $\sum_{i=1}^n\pi^{-\lambda_i}\Nm(a_i)$ lies in the set $F_0/(\mathfrak{o}\cap F_0)$. Notice that $\Nm(\mathfrak{o}^\times)=\mathfrak{o}^\times\cap F_0$, and for all $b\in\mathfrak{o}^\times$ and sequence $(a_1,\ldots,a_n)$ such that $v(\sum_{i=1}^n\pi^{-\lambda_i}\Nm(a_i))=-\lambda_1$, the sequence $(ba_1,\ldots,ba_n)$ satisfies
$$\sum_{i=1}^n\pi^{-\lambda_i}\Nm(ba_i)=\Nm(b)\sum_{i=1}^n\pi^{-\lambda_i}\Nm(a_i)$$
which also has valuation $-\l_1$.
Hence, every element in $F_0/(\mathfrak{o}\cap F_0)$ that has valuation $-\l_1$ (there are $(1-q^{-1})q^{\l_1}$ such elements in total) corresponds to the same number of sequences $(a_1,\ldots,a_n)$. Thus number of such sequences for which $\sum_{i=1}^n\pi^{-\lambda_i}\Nm(a_i)=\pi^{-\lambda_1}$ is given by
$$\frac{q^{2|\lambda|}(1-(-q)^{-m}))(1-q^{-1})}{(1-q^{-1})q^{\lambda_1}}=q^{2|\lambda|-\lambda_1}(1-(-q)^{-m}),$$
as desired.
\end{proof}

Based on \Cref{lem: ways for the first pick}, we have the following theorem, which gives the number of automorphisms for the Hermitian case:

\begin{theorem}\label{thm: number of automorphism_Hermitian}

Let $M_\lambda$ be a Hermitian $\mathfrak{o}$-module of type $\lambda=(\lambda_1,\ldots,\lambda_n)$, where $\lambda_1\ge\ldots\ge\lambda_n\ge 1$. Let $\Aut^H(M_\lambda)$ be the group of module automorphisms of $M_\lambda$ which preserves the pairing on $M_\lambda$. Then
\begin{equation}\#\Aut^H(M_\lambda)=q^{|\lambda|+2n(\lambda)}\prod_{i\ge 1}(-q^{-1};-q^{-1})_{m_i(\lambda)}.\end{equation}
\end{theorem}

\begin{proof}
Since an automorphism may be specified by the images of the standard generators, the number of automorphisms is equal to the number of sequences $x_1,\ldots,x_n\in M_\lambda$ such that
\begin{enumerate}
\item $\langle x_i,x_i\rangle=\pi^{-\lambda_i},\forall 1\le i\le n$, and

\item $\langle x_i,x_j\rangle=0, \forall 1\le i<j\le n$.
\end{enumerate}
We choose $x_1,x_2,\ldots,x_n$ one by one. For the first step, we find an element $x_1\in M_\lambda$ such that $\langle x_1,x_1\rangle=\pi^{-\lambda_1}$. Then, we pick $x_2\in (\mf{o}x_1)^\perp$ such that $\langle x_2,x_2\rangle=\pi^{-\lambda_2}$, and repeat such steps until we reach $x_n$. For every $2\le i\le n$, the element $x_i$ must be chosen in the submodule $(\mf{o}x_1+\cdots+\mf{o}x_{i-1})^\perp$, which is of type $(\lambda_i,\ldots,\lambda_n)$. Applying  \Cref{lem: ways for the first pick}, we write the product of these components, which is 
\begin{align}
\begin{split}
\#\Aut^H(M_\lambda)&=\prod_{i=1}^nq^{\l_i+2\l_{i+1}+\ldots+2\l_n}(1-(-q)^{-\#\{k\ge i\mid\l_j=\l_i\}})\\
&=q^{\sum_i(2i-1)\l_i}\prod_{i\ge 1}(-q^{-1};-q^{-1})_{m_i(\lambda)}\\
&=q^{|\l|+2n(\l)}\prod_{i\ge 1}(-q^{-1};-q^{-1})_{m_i(\lambda)}.
\end{split}
\end{align}
Therefore, we are done.
\end{proof}

Before starting our proof of \Cref{thm: expectation over CL measure}, we introduce an essentially equivalent lemma, which we prove later in this section. 

\begin{lemma}\label[lemma]{lem: sum of skew}
Let $\l,\nu\in\Y$, and $t\in(-1,1) \setminus \{0\}$ be a parameter. Then we have
\begin{equation}\label{eq: sum of skew}
\sum_{\mu\in\Y}\frac{Q_{\l/\mu}(1,t,\ldots;t)P_{\nu/\mu}(t,t^2,\ldots;t)}{Q_{\l}(1,t,\ldots;t)P_{\nu}(t,t^2,\ldots;t)}=t^{-\sum_i\lambda_i'\nu_i'}.
\end{equation}
\end{lemma}

\begin{proof}[Proof of \Cref{thm: expectation over CL measure} from \Cref{lem: sum of skew}]
The desired equality \eqref{eq: expectation over CL measure} is equivalent to the equality
\begin{multline}\label{eq: reduction to common denominator}
P_\nu(t,t^2,\ldots;t)\sum_{\l\in\Y} Q_\l(1,t,\ldots;t)t^{-\sum_i{\l_i'\nu_i'}}P_\l(a_1,a_2,\ldots;t)\\
=P_\nu(a_1,a_2,\ldots,t,t^2,\ldots;t)\Pi_t(a_1,a_2,\ldots;1,t,t^2,\ldots).
\end{multline}
In fact, we have
\begin{align}
\begin{split}
\text{RHS\eqref{eq: reduction to common denominator}}&=\sum_{\mu\in\Y}P_\mu(a_1,a_2,\ldots;t)P_{\nu/\mu}(t,t^2,\ldots;t)\Pi_t(a_1,a_2,\ldots;1,t,t^2,\ldots)\\
&=\sum_{\mu,\l\in\Y}Q_{\l/\mu}(1,t,\ldots;t)P_{\nu/\mu}(t,t^2,\ldots;t)P_\l(a_1,a_2,\ldots;t)\\
&=\sum_{\l\in\Y}Q_{\l}(1,t,\ldots;t)P_{\nu}(t,t^2,\ldots;t)t^{-\sum_i\lambda_i'\nu_i'}P_\l(a_1,a_2,\ldots;t) \\ 
&=\text{LHS\eqref{eq: reduction to common denominator}}.
\end{split}
\end{align}
The second row comes from the skew Cauchy identity, and the third row comes from \eqref{eq: sum of skew}. This ends the proof.
\end{proof}

Now we prove \Cref{lem: sum of skew}, which we used above. The idea of our proof comes from \cite[Proposition 6.3]{nguyen2024universality}, which views the equality \eqref{eq: sum of skew} as counting the maps between abelian $p$-groups.

\begin{proof}[Proof of \Cref{lem: sum of skew}]
Let us start with the case $t=1/p$, where $p$ is a prime number. For all $\l\in\Y$, let  $$G_\l=\bigoplus_{i\ge 1}\Z/p^{\l_i}\Z$$ 
be the abelian $p$-group of type $\l$. The numbers of surjections and injections between abelian $p$-groups are given by
$$\#\Sur(G_\l, G_\mu)=\frac{Q_{\l/\mu}(1,t,\ldots)}{Q_{\l}(1,t,\ldots)P_{\mu}(t,t^2,\ldots)},\quad\#\Inj(G_\mu, G_\nu)=\frac{P_{\nu/\mu}(t,t^2,\ldots)}{Q_{\mu}(1,t,\ldots)P_{\nu}(t,t^2,\ldots)},$$
see \cite[(53)]{nguyen2024universality} for instance. Meanwhile, the number of automorphisms of $G_\mu$ is given by
$$\#\Aut(G_\mu)=\frac{1}{Q_{\mu}(1,t,\ldots)P_{\mu}(t,t^2,\ldots)},$$
see for instance \cite[(52)]{nguyen2024universality}.
Therefore, we have 
\begin{align}\label{eq: from skew sum to hom}
\begin{split}
\text{LHS\eqref{eq: sum of skew}}&=\sum_{\mu\in\Y}\frac{\#\Sur(G_\l, G_\mu)\#\Inj(G_\mu, G_\nu)}{\#\Aut(G_\mu)}\\
&=\sum_{\mu\in\Y}\#\{f\in\Hom(G_\l, G_\nu)\mid f(G_\l)\cong G_\mu\text{ as abelian $p$-group}\}\\
&=\#\Hom(G_\l, G_\nu)\\
&=\prod_{i,j\ge 1}\#\Hom(\Z/p^{\l_i}\Z, \Z/p^{\nu_j}\Z).\\
%&=p^{\sum_{i,j\ge 1}\min(\l_i,\nu_j)}\\
%&=p^{\sum_{r\ge 1} \l_r'\nu_r'}\\
%&=\text{RHS\eqref{eq: sum of skew}},
\end{split}
\end{align}
Notice that $\#\Hom(\Z/p^{\l_i}\Z, \Z/p^{\nu_j}\Z)=p^{\min(\l_i,\nu_j)}$ for all $i,j\ge 1$, and
%where the second-to-last line comes from the equality
$$
\sum_{i,j\ge 1}\min(\l_i,\nu_j)=\sum_{i,j\ge 1}\sum_{r\ge 1}\bbone_{\l_i\ge r,\nu_j\ge r}=\sum_{r\ge 1}\sum_{i,j\ge 1}\bbone_{\l_i\ge r,\nu_j\ge r}=\sum_{r\ge 1} \lambda_r'\nu_r'.
$$
Hence we must have
$$\prod_{i,j\ge 1}\#\Hom(\Z/p^{\l_i}\Z, \Z/p^{\nu_j}\Z)=p^{\sum_{r\ge 1} \l_r'\nu_r'}=\text{RHS\eqref{eq: sum of skew}},$$
and the equality \eqref{eq: sum of skew} must hold when $t=1/p$. Since there exist infinitely many prime numbers, and both sides of \eqref{eq: sum of skew} belong to $\Q(t)$, the same equality must hold when we view $t$ as a parameter.
\end{proof}

One can also prove \Cref{lem: sum of skew} by direct computation and induction; see \Cref{sec:appendix}.

\begin{proof}[Proof of \Cref{thm: expectation of hom}] 
For the no-pairing case, when the usual module $M$ has type $\l$, we have $\#\Hom(M,M_\nu)=q^{\sum_i\l_i'\nu_i'}$, as showed in \eqref{eq: from skew sum to hom}. Applying \Cref{thm: expectation over CL measure} and setting $t=q^{-1}$ brings the result. 

For the alternating case, when the alternating module $M$ has type $\l$, we have $\#\Hom(M,M_\nu)=q^{\sum_i2\l_i'\nu_i'}=(q^2)^{\sum_i\l_i'\nu_i'}$. Applying \Cref{thm: expectation over CL measure} and setting $t=q^{-2}$ brings the result. 

Likewise, for the Hermitian case, when the Hermitian module $M$ has type $\l$, we have $\#\Hom(M,M_\nu)=(q^2)^{\sum_i\l_i'\nu_i'}=(-q)^{\sum_i\l_i'(\nu^2)_i'}$. Applying \Cref{thm: expectation over CL measure} and setting $t=-q^{-1}$ brings the result.

\end{proof}

\begin{rmk} The $u$-measure on groups without pairings corresponds to the case $a_1=t^{1+u},a_2=t^{2+u},\ldots$ in the notation of \Cref{eq:hl_measure_intro}, where $t=1/p$. In this case, the Hall-Littlewood formulas for the Hom-moments in \eqref{eq: expectation over CL measure} have another form as a series expansion in subgroup counts. Namely,
\begin{align}
\begin{split}
\text{RHS\eqref{eq: expectation over CL measure}}&=\frac{P_\nu(t^u,t^{1+u},\ldots,1,t,\ldots)}{P_\nu(1,t,\ldots)}\\
&=\sum_{\mu\in\Y}\frac{P_\mu(t^u,t^{1+u},\ldots)}{P_\nu(1,t,\ldots)}P_{\nu/\mu}(1,t,\ldots)
\\
&=\sum_{\mu\in\Y}\frac{P_\mu(1,t,\ldots)}{P_\nu(1,t,\ldots)}P_{\nu/\mu}(1,t,\ldots)t^{|\mu|u} \\ 
&=\sum_{\mu\in\Y}C_{\nu/\mu}(p)p^{-|\mu|u},
\end{split}
\end{align}
where $C_{\nu/\mu}(p)=\frac{P_\mu(1,t,\ldots)}{P_\nu(1,t,\ldots)}P_{\nu/\mu}(1,t,\ldots)$ is the number of subgroups in $G_\nu$ of type $\mu$, as discussed in the remark on \cite[page 12]{delaunay2014p}. This was proven in \cite[Theorem 12]{delaunay2014p}, by a different and more involved method.
\end{rmk}

\appendix

\section{A direct proof of {\Cref{lem: sum of skew}}}\label{sec:appendix}

In this appendix, we prove \Cref{lem: sum of skew} by direct computation. For all $n\in\Z$ and $m\ge 0$, let 
$$\begin{bmatrix} n\\m\end{bmatrix}_t=\frac{(1-t^{n-m+1})\cdots(1-t^n)}{(1-t)\ldots(1-t^m)}$$
be the $q$-binomials. Let us start from the following lemma.

\begin{lemma}\label[lemma]{lem: length one case}
Let $n,\l_1,\nu_1\ge 0$ be integers. Denote 
\begin{equation}\label{eq: def_f}
    f(n;\l_1,\nu_1)=\sum_{\mu_1}t^{-(\mu_1+n)(\l_1+\nu_1-\mu_1+n)}(t^{1+\l_1-\mu_1};t)_{\mu_1}\begin{bmatrix} \nu_1\\\mu_1\end{bmatrix}_t.
\end{equation}
Then we have $f(n;\l_1,\nu_1) = t^{-(\l_1+n)(\nu_1+n)}$.
\end{lemma}

\begin{proof}
We prove the desired equality $f(n;\l_1,\nu_1) = t^{-(\l_1+n)(\nu_1+n)}$ by induction on $\nu_1$. It is clear the lemma holds when $\nu_1=0$. Suppose the equation already holds for some $\nu_1$. Then we have
\begin{align}
\begin{split}
f(n;\lambda_1,\nu_1+1)&=\sum_{\mu_1}t^{-(\mu_1+n)(\l_1+\nu_1-\mu_1+n+1)}(t^{1+\l_1-\mu_1};t)_{\mu_1}\begin{bmatrix} \nu_1+1\\\mu_1\end{bmatrix}_t\\
&=\sum_{\mu_1}t^{-(\mu_1+n)(\l_1+\nu_1-\mu_1+n+1)}(t^{1+\l_1-\mu_1};t)_{\mu_1}(t^{\mu_1}\begin{bmatrix} \nu_1\\\mu_1\end{bmatrix}_t+\begin{bmatrix} \nu_1\\\mu_1-1\end{bmatrix}_t)\\
&=t^{-n}f(n;\l_1,\nu_1)+\sum_{\mu_1}t^{-(\mu_1+n)(\l_1+\nu_1-\mu_1+n+1)}(t^{1+\l_1-\mu_1};t)_{\mu_1}\begin{bmatrix} \nu_1\\\mu_1-1\end{bmatrix}_t\\
&=t^{-n}f(n;\l_1,\nu_1)+\sum_{\mu_1}t^{-(\mu_1+n+1)(\l_1-1+\nu_1-\mu_1+n+1)}(1-t^{\l_1})(t^{\l_1-\mu_1};t)_{\mu_1}\begin{bmatrix} \nu_1\\\mu_1\end{bmatrix}_t\\
&=t^{-n}f(n;\l_1,\nu_1)+(1-t^{\l_1})f(n+1;\l_1-1,\nu_1)\\
&=t^{-n-(\l_1+n)(\nu_1+n)}+(1-t^{\l_1})t^{-(\l_1+n)(\nu_1+n+1)}\\
&=t^{-(\l_1+n)(\nu_1+n+1)},
\end{split}
\end{align}
where the fourth row is deduced from the third row by replacing $\mu_1$ with $\mu_1+1$. Therefore, the lemma also holds for $\lambda_1,\nu_1+1$. This ends the proof. 
\end{proof}

Based on the above preparation, we now turn back to the proof of \Cref{lem: sum of skew}.

\begin{proof} 
By the principal specialization formulas \eqref{eq: principal specialization formula for P}, \eqref{eq: principal specialization formula for Q}, this is equivalent to
$$\sum_{\mu\in\Y}t^{\sum_i\left(\begin{pmatrix} \lambda_i'-\mu_i'\\ 2\end{pmatrix}+\begin{pmatrix} \nu_i'-\mu_i'\\ 2\end{pmatrix}-\begin{pmatrix} \lambda_i'\\ 2\end{pmatrix}-\begin{pmatrix} \nu_i'\\ 2\end{pmatrix}-\mu_i'\right)}\prod_{i\ge 1}\frac{(t^{1+\l_i'-\mu_i'};t)_{m_i(\mu)}(t^{1+\nu_i'-\mu_i'};t)_{m_i(\mu)}}{(t;t)_{m_i(\mu)}}=t^{-\sum_i\lambda_i'\nu_i'}.$$
Notice the fact that
$$\begin{pmatrix} \lambda_i'-\mu_i'\\ 2\end{pmatrix}+\begin{pmatrix} \nu_i'-\mu_i'\\ 2\end{pmatrix}-\begin{pmatrix} \lambda_i'\\ 2\end{pmatrix}-\begin{pmatrix} \nu_i'\\ 2\end{pmatrix}-\mu_i'=-\mu_i'(\lambda_i'+\nu_i'-\mu_i').$$ 
Hence we only need to prove conjugate version, which has the form
\begin{equation}\label{eq: arbitrary length case}
\sum_{\mu\in\Y}t^{-\sum_{i\ge 1} \mu_i(\lambda_i+\nu_i-\mu_i)}\prod_{i\ge 1}\frac{(t^{1+\l_i-\mu_i};t)_{\mu_i-\mu_{i+1}}(t^{1+\nu_i-\mu_i};t)_{\mu_i-\mu_{i+1}}}{(t;t)_{\mu_i-\mu_{i+1}}}=t^{-\sum_i\lambda_i\nu_i}.
\end{equation}

Notice that 
\begin{align}
\begin{split}
\text{LHS\eqref{eq: arbitrary length case}}&=\sum_{\mu_2\ge\mu_3\ge\ldots}t^{-\sum_{i\ge 2} \mu_i(\lambda_i+\nu_i-\mu_i)}\prod_{i\ge 2}(t^{1+\l_i-\mu_i};t)_{\mu_i-\mu_{i+1}}\begin{bmatrix} \nu_i-\mu_{i+1} \\ \mu_i-\mu_{i+1}\end{bmatrix}_t\cdot\\
&\sum_{\mu_1\ge \mu_2}t^{-\mu_1(\lambda_1+\nu_1-\mu_1)}(t^{1+\l_1-\mu_1};t)_{\mu_1-\mu_2}\begin{bmatrix} \nu_1-\mu_2 \\ \mu_1-\mu_2\end{bmatrix}_t\\
&=\sum_{\mu_2\ge\mu_3\ge\ldots}f(\mu_2;\l_1-\mu_2,\nu_1-\mu_2)t^{-\sum_{i\ge 2} \mu_i(\lambda_i+\nu_i-\mu_i)}\prod_{i\ge 2}(t^{1+\l_i-\mu_i};t)_{\mu_i-\mu_{i+1}}\begin{bmatrix} \nu_i-\mu_{i+1} \\ \mu_i-\mu_{i+1}\end{bmatrix}_t\\
&=t^{-\l_1\nu_1}\sum_{\mu_2\ge\mu_3\ge\ldots}t^{-\sum_{i\ge 2} \mu_i(\lambda_i+\nu_i-\mu_i)}\prod_{i\ge 2}(t^{1+\l_i-\mu_i};t)_{\mu_i-\mu_{i+1}}\begin{bmatrix} \nu_i-\mu_{i+1} \\ \mu_i-\mu_{i+1}\end{bmatrix}_t,\\
\end{split}
\end{align}
where the second equality is obtained by replaing the variable $\mu_1-\mu_2$ with $\mu_2$, and the third equality is obtained by \Cref{lem: length one case}. Iterating this method in a similar way and we have 
\begin{align}
\begin{split}
%\text{LHS\eqref{eq: arbitrary length case}}&=\sum_{\mu_2,\mu_3,\ldots}f(\mu_2;\l_1-\mu_2,\nu_1-\mu_2)t^{-\sum_{i\ge 2} \mu_i(\lambda_i+\nu_i-\mu_i)}\prod_{i\ge 2}(t^{1+\l_i-\mu_i};t)_{\mu_i-\mu_{i+1}}\begin{bmatrix} \nu_i-\mu_{i+1} \\ \mu_i-\mu_{i+1}\end{bmatrix}_t\\
\text{LHS\eqref{eq: arbitrary length case}}&=t^{-\l_1\nu_1}\sum_{\mu_2,\mu_3,\ldots}t^{-\sum_{i\ge 2} \mu_i(\lambda_i+\nu_i-\mu_i)}\prod_{i\ge 2}(t^{1+\l_i-\mu_i};t)_{\mu_i-\mu_{i+1}}\begin{bmatrix} \nu_i-\mu_{i+1} \\ \mu_i-\mu_{i+1}\end{bmatrix}_t\\
&=t^{-\l_1\nu_1-\l_2\nu_2}\sum_{\mu_3,\mu_4,\ldots}t^{-\sum_{i\ge 3} \mu_i(\lambda_i+\nu_i-\mu_i)}\prod_{i\ge 3}(t^{1+\l_i-\mu_i};t)_{\mu_i-\mu_{i+1}}\begin{bmatrix} \nu_i-\mu_{i+1} \\ \mu_i-\mu_{i+1}\end{bmatrix}_t\\
&=\ldots\\
&=t^{-\sum_i\l_i\nu_i}\\
&=\text{RHS\eqref{eq: arbitrary length case}}.
\end{split}
\end{align}
Only finitely many iterations are necessary because both $\l$ and $\nu$ are of finite length. This ends the proof. 
\end{proof}

\textbf{Acknowledgements.} We thank Cesar Cuenca, Jason Fulman, Chao Li, and Grigori Olshanski for helpful conversations, and the anonymous referee for many helpful suggestions. JS was partially supported by NSF grant DMS-2246576 and Simons Investigator
grant 929852.

\end{document}